\documentclass[a4paper,10pt]{article}
\usepackage{a4,amsmath,amssymb,enumerate,stmaryrd,verbatim}
\usepackage[dvips,bookmarksnumbered=true]{hyperref}
\renewcommand{\epsilon}{\varepsilon}
\renewcommand{\phi}{\varphi}
\renewcommand{\rho}{\varrho}
\def\0{\relax\ifmmode\emptyset\else$\emptyset$\fi}
\def\ge{\geqslant} 
\def\le{\leqslant} 
\def\models{\vDash}

\newcommand{\mbb}{\mathbb}
\newcommand{\mc}{\mathcal}
\newcommand{\mf}{\mathfrak}
\newcommand{\mr}{\mathrm}

\newcommand{\tb}{\textbf}

\newcommand{\Aut}{\mathrm{Aut}}
\newcommand{\elem}{\preccurlyeq}
\newcommand{\falsum}{\bot}

\newcommand{\into}{\hookrightarrow} 
\newcommand{\isom}{\cong} 
\newcommand{\liff}{\leftrightarrow}

\newcommand{\onto}{\twoheadrightarrow} 
\newcommand{\verum}{\top}

\newenvironment{List}[3][$\bullet$]{\begin{list}{#1}{\itemsep#3ex\setlength{\topsep}{0.5\itemsep}\parsep0ex\labelsep1ex\settowidth{\labelwidth}{#1}\setlength{\leftmargin}{\labelwidth}\addtolength{\leftmargin}{#2ex}\addtolength{\leftmargin}{\labelsep}}}{\end{list}} 

\long\def\footnotestar#1{\begingroup%
\def\thefootnote{\fnsymbol{footnote}}\footnote[0]{#1}\endgroup}

\newcommand{\longonto}{-\!\!\!\onto}
\newcommand{\loongonto}{-\!\!\!-\!\!\!\onto}
\newcommand{\down}{\downarrow}
\newcommand{\up}{\uparrow}
\newcommand{\erz}[1]{{#1}_\down}
\newcommand{\zre}[1]{{#1}^{\mbox{\protect\raisebox{1.3pt}{$\scriptstyle\up$}}\complement}}
\newcommand{\Erz}[1]{\{#1\}_\down}
\newcommand{\Zre}[1]{\{#1\}^{\mbox{\protect\raisebox{1.3pt}{$\scriptstyle\up$}}\complement}}
\newcommand{\cpl}{^\complement}
\newcommand{\dual}{^{*}} 
\newcommand{\rin}{\mbox{\raisebox{0.7pt}{$\scriptstyle\in$}}}
\newcommand{\Spec}{\mathrm{Spec}}
\newcommand{\supp}{\mathrm{supp}}
\newcommand{\suppmin}{\mathrm{supp}^{\mathrm{min}}}

\newcommand{\val}{\mr{val}}
\newcommand{\Val}{\mr{Val}}
\newcommand{\wh}{\widehat} 
\newcommand{\Fr}[2][]{\mc F^{#1}_{#2}}
\newcommand{\Fco}[1]{\widehat{\mc F^{}_{#1}}} 
\newcommand{\Fkl}[1]{\smash[t]{\overset{\text{\raisebox{-1.5pt}[1.5pt][0pt]{$\smile$}}}{\mc F^{}_{#1}}}}
\newcommand{\IPL}[1]{$\mathrm{IPL}_{#1}$}

\newcommand{\hto}{\to}
\newcommand{\hiff}{\leftrightarrow}
\newcommand{\hle}{\sqsubseteq}
\newcommand{\hge}{\sqsupseteq}
\newcommand{\hleneq}{\sqsubset}

\newcommand{\hcap}{\sqcap}
\newcommand{\hcup}{\sqcup}
\newcommand{\bighcap}{\bigsqcap}
\newcommand{\bighcup}{\bigsqcup}
\newcommand{\Hcap}{\bigsqcap}
\newcommand{\Hcup}{\bigsqcup}

\newtheorem{Theorem}{Theorem}[section]
\newtheorem{Cor}[Theorem]{Corollary}

\newtheorem{EX}[Theorem]{Example}
\newtheorem{Fact}[Theorem]{Fact}
\newtheorem{Lemma}[Theorem]{Lemma}
\newtheorem{Prop}[Theorem]{Proposition}
\newtheorem{QU}[Theorem]{Question}
\newtheorem{REM}[Theorem]{Remark}
\newtheorem{Problem}{Problem}
\newenvironment{Rem}{\begin{REM}\normalfont}{\end{REM}}
\newenvironment{Ex}{\begin{EX}\normalfont}{\end{EX}}
\newenvironment{Qu}{\begin{QU}\normalfont}{\end{QU}}

\newcommand{\PROOF}[1][]{\textsc{Proof#1:} }
\newcommand{\END}{\hfill$\square$\medskip}

\begin{document}

\title{
  On Bellissima's construction of the finitely generated free Heyting
  algebras, and beyond
  } 
\author{
  \begin{minipage}[t]{5.5cm} 
  Luck Darni\`ere 
   \thanks{The first author would like to thank the Universit\"at
     Freiburg for inviting him in July 2008.} \\ 
   \begin{small} 
     D\'epartement de math\'ematiques \\[-1ex] 
     Facult\'e des sciences \\[-1ex] 
     Universit\'e d'Angers, France 
   \end{small}
  \end{minipage}
  \begin{minipage}[t]{5.5cm} 
  Markus Junker\addtocounter{footnote}{5}
   \thanks{The second author would like to thank the Universit\'e 
     d'Angers for supporting him as an invited professor in march
     2005, when main parts of this work were done.} \\
   \begin{small}
     Mathematisches Institut \\[-1ex] 
     Abteilung f\"ur mathematische Logik \\[-1ex] 
     Universit\"at Freiburg, Deutschland 
   \end{small}
  \end{minipage}
  }

\date{December 8, 2008}

\maketitle

\footnotestar{{\bf MSC 2000:} 06D20, 03C64, 06B23, 06B30, 08B20}

\section{Introduction}

Heyting algebras are a generalisation of Boolean algebras; the most
typical example is the lattice of open sets of a topological space. 
Heyting algebras play the same r\^ole for intuitionistic logic as
Boolean algebras for classical logic. They are special distributive
lattices, and they form a variety. They are mainly studied by universal
algebraists and by logicians, hardly by model theorists. In contrast 
to Boolean algebras, finitely generated free Heyting algebras are
infinite, as was shown in the first article on Heyting algebras by
McKinsey and Tarski in the 1940s. 
For one generator, the free Heyting algebra is well understood, but
from two generators on, the structure remains mysterious, though many
properties are known. With the help of recursively described Kripke
models, Bellissima has given a representation of the finitely
generated free Heyting algebras $\Fr{n}$ as sub-algebras of
completions $\Fco{n}$ of them. Essentially the same construction is
due independently to Grigolia. Our paper offers a concise and readable
account of Bellissima's construction and analyses the situation closer. 

Our interest in Heyting algebras comes from model theory and geometry.
Our initial questions concerned axiomatisability and decidability of
structures like the lattice of Zariski closed subsets of $K^n$ for
fields $K$, which led us rapidly to questions about Heyting algebras. 
One of the problems with Heyting algebras is that they touch many
subjects: logic, topology, lattice theory, universal algebra, category
theory, computer science. Therefore there are many different
approaches and special languages, which often produce papers that are
hard to read for non-insiders. An advantage of our article should be
clear proofs and the use mainly of standard mathematical
terminology. There is a bit of logic that one might skip if one
believes in Bellissima's theorem; and there are basic model theoretic
notions involved in the section about model theoretic results. 

Our paper is organised as follows:
Section~\ref{Sec-basic} introduces the definitions, basic properties,
and reference examples. Section~\ref{Sec-Bellissima} contains an
account of Bellissima's construction, a short proof, and results from
his article that we are not going to prove. Section~\ref{Sec-csq}
analyses the Heyting algebra constructed by Bellissima: we show that
it is the profinite completion of the finitely generated free Heyting
algebra as well as the metric completion for a naturally defined metric. 
In Section~\ref{Sec-vier}, we reconstruct the Kripke model as
the principal ideal spectrum and prove that the Zariski topology on
this spectrum is induced by the partial ordering. Section~\ref{Sec-mt}
shows the Kripke model to be first order interpretable, from which
several model theoretic and algebraic properties for dense
sub-algebras of $\Fco{n}$ follow. For example, we show the set of
generators to be $\0$-definable, and we determine automorphism groups.
We solve questions of elementary equivalence, e.g.\ we prove that no 
proper sub-algebra of $\Fr{n}$ is elementarily equivalent to $\Fr{n}$,
and that $\Fr{n}$ is an elementary substructure of $\Fco{n}$ iff both
algebras are elementarily equivalent. And we settle some questions about 
irreducible elements: $\Fr{n}$ contains the same meet-irreducible
elements as $\Fco{n}$; we characterise the join-irreducible elements,
we show that there are continuum many of them in $\Fco{n}$ and that
the join-irreducibles of $\Fr{n}$ remain join-irreducible in $\Fco{n}$.
Finally, Section~\ref{Sec-letzt} collects open problems and
miscellaneous considerations.  

Some of the properties we isolated were known before, some were
published after we started this work in 2004. We added references
where we were able to do so, and apologise for everything we have
overlooked. Though not all the results are new, the proofs might be, 
and the way of looking at the problem is hopefully interesting. 

This paper is closely related to \cite{DJ2}; both complement each
other. When we started to study finitely generated Heyting algebras,
we did it in two ways: on the one hand by analysing Bellissima's
construction, on the other hand by analysing the notion of dimension
and codimension in dual Heyting algebras. Many insights were obtained
by both approaches, but some features are proper to the free Heyting
algebras, others hold for a much wider class than just the finitely
generated Heyting algebras. Therefore, we decided to write two papers:
this one, which collects results that follow more or less directly
from Bellissima's construction, and the paper \cite{DJ2}, which
analyses the structure of Heyting algebras from a more geometric point
of view.

\subsubsection*{Acknowledgements}

We would like to thank Guram Bezhanishvili for many helpful
remarks.

\section{Basic facts about Heyting algebras} \label{Sec-basic}

\subsection{Definitions and notations}

Under a \emph{lattice} we will always understand a distributive
lattice with maximum $1$ and minimum $0$, join (union) and meet 
(intersection) being denoted by $\hcup$ and $\hcap$ respectively. 
A \emph{Heyting algebra} is a lattice where for every $a,b$ there
exists an element  
$$a \hto b \ :=\ \max \big\{ x \bigm| x \hcap a = b\hcap a \big\}.$$ 
This is expressible as a
universal theory, the theory $T_{HA}$ of Heyting algebras, in the
language $\mc L_{HA} = \{0,1,\hcap,\hcup,\hto\}$ with constant
symbols $0,1$ and binary function symbols $\hcap,\hcup,\hto$. 

Note that everything is definable from the partial ordering 
$$ a \hle b \ :\iff\  a = a\hcap b 
           \  \iff\   b = a\hcup b 
           \  \iff\   a\hto b = 1$$

In a poset $(X,\le)$, we call $y$ \emph{a successor} of $x$ if $x<y$ 
and there is no $z$ with $x<z<y$. Analogously for \emph{predecessor}. 

If $\Lambda = (\Lambda,0,1,\hcap,\hcup,\hle)$ is a lattice, then the
dual lattice $\Lambda\dual := (\Lambda,1,0,\hcup,\hcap,\hge)$ is the
lattice of the reversed ordering. Thus, in the dual of a Heyting
algebra, for all $a,b$ there exists a smallest element $b-a$ with the
property $a \,\hcup\, (b-a) = a \,\hcup\, b$. Dual Hey\-ting algebras 
are the older siblings of Heyting algebras; they were born
\emph{Brouwerian algebras} in \cite{McTarski}, they appear under the
name \emph{topologically complemented lattices} in \cite{D-preprint},
and are often called \emph{co-Heyting algebras} nowadays. A lattice is
a \emph{bi-Heyting algebra} (or \emph{double Brouwerian algebra} in
\cite{McTarski}) if itself and its dual are Heyting algebras. 

We define $a \hiff b$ as $(a \hto b) \hcap (b \hto a)$, which is dual
to the ``symmetric difference'' $a \vartriangle b := (a-b) \hcup (b-a)$.

\subsection{Examples} \label{Ex}

\begin{List}12
\item 
  The open sets of a topological space $X$ form a complete%
 \footnote{\textit{I.e.} complete as a lattice. Note that in general only one
    of the infinite distributive laws holds in topological Heyting algebras, namely 
    $a \hcap \bighcup_{i\in I} b_i = \bighcup_{i\in I} (a \hcap b_i)$.}
  Heyting algebra $\mc O(X)$ with the operations suggested by the
  notations, i.e. $a \hcap b = a \cap b$, $a \hcup b = a \cup b$, $a
  \hle b \iff a \subseteq b$, and 
  $$a \hto b\ =\ {\overline{a \setminus b}}\cpl\ =\ (a\cpl \cup b)^\circ.$$ 
  (Here, ${}\cpl$ denotes the complement, $\overline{\phantom{M}}$ topological
  closure and ${ }^\circ$ the interior.)\\ 
  We call such an algebra a \emph{topological Heyting algebra}. 

\item 
  If $(X, \le)$ is a partial ordering, then the increasing sets form a 
  topology and hence a Heyting algebra $\mc O^\up(X,\le)$, and the
  decreasing sets, which are the closed sets of $\mc O^\up(X,\le)$,
  form a topology and Heyting algebra $\mc O_\down(X,\le)$. Hence such
  an algebra is bi-Heyting, and both infinite distributive laws hold. 

\item 
  The propositional formulae in $\kappa$ propositional variables%
 \footnote{For these ``intuitionistic formulae'', the system of
    connectives $\{\falsum,\land,\lor,\to\}$ is used, and the following 
    abbreviations: $\verum := \falsum\to\falsum$, $\lnot A := A \to \falsum$, 
    $A \liff B := (A \to B)\land(B\to A)$.},
  up to equivalence in the intuitionistic propositional calculus, form
  a Heyting algebra \IPL{\kappa}. It is freely generated by the
  (equivalence classes of the) propositional variables, hence
  isomorphic to the \emph{free Heyting algebra $\Fr{\kappa}$ over
    $\kappa$ generators}. 

 \item 
  If $\pi: H \to H'$ is a non-trivial epimorphism of Heyting algebras,
  then the \emph{kernel} $\pi^{-1}(1)$ is a filter.
  Conversely, if $\Phi$ is an arbitrary filter in $H$, then
  $\equiv_\Phi$ defined by $x \equiv_\Phi y :\iff x \hiff y \in \Phi$
  is a congruence relation such that $\Phi = \pi_\Phi^{-1}(1)$ for the 
  canonical epimorphism $\pi_\Phi: H \to H/{\equiv_\Phi}$. 
  In particular, $\equiv_{\{1\}}$ is equality. 

\end{List}

\section{Bellissima's construction} \label{Sec-Bellissima}

Bellissima in \cite{B} has constructed an embedding of the free
Heyting algebra $\Fr{n}$ into the Heyting algebra $\mc O_\down(\mf K_n)$ 
for a ``generic'' Kripke model $\mf K_n$ of intuitionistic
propositional logic in $n$ propositional variables $P_1,\dots,P_n$. 
We fix $n>0$ and this fragment \IPL{n} of intuitionistic logic, and we
will give a short and concise account of Bellissima's construction
and proof. 

We start with some terminology: 
we identify the set $\Val_n$ of all valuations (assignments) of the
propositional variables with the power set of $\{P_1,\dots,P_n\}$,
namely a valuation with the set of variables to which it assigns ``true''. 
A \emph{Kripke model} $\mf K = (K,\le,\val)$ for \IPL{n} consists of a
reflexive partial order%
 \footnote{Note that the order is reversed with respect to the usual
  approach to Kripke models. This is for the sake of an easy
  description and to be in coherence with the order of the Heyting
  algebra, see Remark~\ref{Rem-suppKripke}, and the order on the
  spectrum, compare with Fact~\ref{F-spectrum}.} 
$\le$ on a set $K$ and a function $\val:K \to \Val_n$
satisfying the following monotonicity condition: If $P_i \in \val(w)$
for $w \in K$, written $w\models P_i$, and if $w' \le w$, then $w'\models P_i$.
Validity of formulae at a point $w$ is then defined by induction in
the classical way for the connectives $\falsum$, $\lor$ and $\land$,
and for $\phi\to\chi$ by the condition: $w'\models\phi \Rightarrow w'
\models \chi$ for all points $w' \in K$ with $w'\le w$. By induction,
validity of all intuitionistic formulae obeys the monotonicity condition. 
It follows from this definition that the map 
$$\phi\ \mapsto\ [\![\phi]\!] \,:=\, \{ w \in K \mid w\models \phi\}$$ 
induces a homomorphism of Heyting algebras from \IPL{n} to 
$\mc O_\down(K,\le)$. The kernel of this morphism is the 
\emph{theory of the model}: all formulae valid at every point of the model. 
The \emph{theory of a point} $w$ consists of all formulae valid at $w$. 
(See e.g.\ \cite{Fitting} for more details.) 

A Kripke model is \emph{reduced} if any two points differ either by
their valuations or by some (third) point below. Precisely: there are
no two distinct points $w_1,w_2$ with the same valuation and (1) such
that $w \le w_1 \iff w \le w_2$ for all $w$ or (2) such that $w_1$ is
the unique predecessor of $w_2$. (One can show that a finite model is
reduced if and only if two distinct points have distinct theories.)
One can reduce a finite model by applying the following two operations: 
\begin{List}20
\item 
  identify points with same valuation and same points below;
\item
  delete a point with only one predecessor, if both 
  carry the same valuation; 
\end{List} 
and one can check by induction that these reductions do not change the 
theory of the model. (This is well known in modal logic: it is a
special case of a bisimulation, see \cite{BMV}.) Therefore it follows:  

\begin{Fact}[Lemma 2.3 of \cite{B}] \label{F-reduced} 
  For every finite model of\/ \IPL{n},
  there is a reduced finite model with the same theory.%
\end{Fact}

\subsection{The construction of the generic Kripke model $\mf K_n$}

The idea of the construction of $\mf K_n$ is to ensure that all finite
reduced models embed as an initial segment. $(K_n,\le)$ will be a 
well-founded partial ordering of rank $\omega$ and $\mf K_n$ will be
an increasing union of Kripke models $\mf K_n^d = (K_n^d,\le,\val)$. 
We define $\mf K_n^d$ by induction on $d$ as follows (cf.\ Figure~\ref{f1}):

\unitlength1cm
\begin{figure}[ht] 
\begin{center}
\begin{picture}(12,5.5)
\put(2.5,0){\line(-1,5){1}}
\put(8.5,0){\line(1,5){1}}
\put(2.5,0){\line(1,0){6}}
\put(2.4,0.5){\line(1,0){6.2}}
\put(2.3,1){\line(1,0){6.4}}
\put(2.2,1.5){\line(1,0){6.6}}
\put(5.5,1.8){$\vdots$}
\put(2,2.5){\line(1,0){7}}
\put(1.9,3){\line(1,0){7.2}}
\put(1.8,3.5){\line(1,0){7.4}}
\put(5.5,3.8){$\vdots$}
\put(9.5,0.1){$K_n^0$}
\put(9.6,0.6){$K_n^1\setminus K_n^0$}
\put(9.7,1.1){$K_n^2 \setminus K_n^1$}
\multiput(10.2,1.85)(0.02,0.1){3}{.}
\put(10,2.6){$K_n^d\setminus K_n^{d-1}$}
\put(10.1,3.1){$K_n^{d+1} \setminus K_n^d$}
\multiput(10.6,3.85)(0.02,0.1){3}{.}
\put(0,1.3){$K_n^d$}
\put(1.5,0.25){\oval(0.5,0.5)[bl]}
\put(1,1.25){\oval(0.5,0.5)[tr]}
\put(1,1.75){\oval(0.5,0.5)[br]}
\put(1.5,2.75){\oval(0.5,0.5)[tl]}
\put(1.25,0.25){\line(0,1){1}}
\put(1.25,1.75){\line(0,1){1}}
\put(3,0.2){$w_{\0,\0}$}
\put(2.8,0.3){\circle*{0.1}}
\put(3.8,0.2){$\cdots$}
\put(7.5,3.2){$w_{\beta,Y}$}
\put(7.3,3.3){\circle*{0.1}}
\thicklines
\put(4.5,0){\line(1,2){0.5}}
\put(4.52,0){\line(1,2){0.5}}
\put(5,1.02){\line(1,0){1}}
\put(6,1){\line(1,2){1}}
\put(6.02,1){\line(1,2){1}}
\put(7,3.02){\line(1,0){0.5}}
\put(7.48,3){\line(1,-3){1}}
\put(7.5,3){\line(1,-3){1}}
\put(7,1.8){$Y$}
\end{picture} 
\end{center}
\caption{The construction of $\mf K_n$\label{f1}}
\end{figure}

\begin{List}21
\item
  We let $K_n^{-1}= \emptyset$. Then $K_n^d\setminus K_n^{d-1}$ consists of
  all possible elements $w_{\beta,Y}$ such that: 
  \begin{List}[$\circ$]20
  \item 
    $Y$ is a decreasing set in $K_n^{d-1}$ and $Y \not\subseteq K_n^{d-2}$ \\ 
    (for $d=0$ the last condition is empty, therefore $Y=\emptyset$);
  \item 
    $\beta$ is a valuation in $\Val_n$ such that $\beta \subseteq
    \val(w')$ for all points $w' \in Y$;
  \item 
    if $Y$ is the decreasing set generated by an element
    $w_{\beta',Y'}$, then $\beta \neq \beta'$.
  \end{List}
\item
  The valuation of $w_{\beta,Y}$ is defined to be $\beta$. 
\item 
  The partial ordering on $K_n^{d-1}$ is extended to $K_n^d$ by
  $$w \le w_{\beta,Y} :\iff (w \in Y \text{ or } w=w_{\beta,Y}).$$
\end{List}

In particular, one sees that by construction every $K_n^d$ is finite,
$K_n^d$ is an initial part of $K_n^{d+1}$, and $K_n^d$ is the set of
points of $K_n$ of foundation rank $\le d$. One can check that $\mf K_n^d$ 
is the maximal reduced Kripke model of foundation rank $d$ for \IPL{n}.

\begin{Theorem}[Bellissima in \cite{B}] \label{Thm-B1} 
The map $\phi \mapsto [\![\phi]\!] = \{ w \in K_n \mid w\models \phi\}$ 
induces an embedding of \IPL{n}, and hence of the free Heyting algebra
$\Fr{n}$ with a fixed enumeration of $n$ free generators, into the
Heyting algebra $$\Fco{n}\ :=\ \mc O_\down(K_n,\le).$$ 
\end{Theorem} 

This map identifies a set of $n$ free generators of $\Fr{n}$ with the
propositional variables $P_1,\dots,P_n$ that were used in the
construction of $\mf K_n$. We will see in Corollary~\ref{Cor-freegen}
that there is only one set of free generators of $\Fr{n}$, therefore
the embedding is unique up to the action of $\mr{Sym}(n)$ on the free
generators.

\medskip 
\PROOF[\footnote{The proof is essentially Bellissima's: we have
  simplified notations, separated the general facts from the special
  situation and left out some detailed elaborations, e.g. a proof of
  Fact~\ref{F-reduced}.}]
We have to show that the homomorphism is injective, which amounts to
show that the theory of $\mf K_n$ consists exactly of all intuitionistic 
tautologies. Each finite reduced Kripke model embeds by a
straightforward induction onto an initial segment of $(\mf K_n,\le)$. 
Because validity of formulae is preserved under ``going down'' along
$\le$, the theory of $\mf K_n$ is contained in that of all finite
reduced models. On the other hand, as intuitionistic logic has the
finite model property, a non-tautology is already false in some finite
model, hence also in some finite reduced model. Thus the theory of
$\mf K_n$ consists exactly of the intuitionistic tautologies. 
\END

For Grigolia's version of this construction see e.g. \cite{Grigolia-free} 
or the account in \cite{Nick} which offers a wider context.

\medskip 
From now on, we will identify $\Fr{n}$ with its image in $\Fco{n}$.
Thus the free generators become $[\![P_1]\!],\dots,[\![P_n]\!]$, and
the operations can be computed in the topological Heyting algebra
$\Fco{n}$ as indicated in Example~\ref{Ex}. Moreover, we speak of
\emph{finite elements} of $\Fr{n}$ or $\Fco{n}$ meaning elements which
are finite subsets of $K_n$.

\subsection{Results from \cite{B}}

In this section we collect all results from \cite{B} that we are going
to use. 

We call an element $a$ in a lattice \emph{$\hcup$-irreducible} if it
is different from $0$ and can not be written as a union $b_1 \hcup
b_2$ with $b_i \neq a$. It is \emph{completely $\hcup$-irreducible},
or \emph{$\Hcup$-irreducible} for short, if it can not be written as
any proper union of other elements (possibly infinite, possibly empty), 
thus if and only if it has a unique predecessor $a^- := \bighcup\{x \mid x<a\}$. 
The dual notions apply for $\hcap$. In particular, a $\hcap$-irreducible 
element is by convention different from $1$, and a $\Hcap$-irreducible
has a unique successor $a^+ := \bighcap \{ x \mid a < x\}$. 


For $X \subseteq K_n$, we let $\erz X$ be the $\le$-decreasing set
and $X^\up$ the $\le$-increasing set generated by $X$.%
\footnote{Because we are working here with the reversed order, the
  arrows are the other way round compared to Bellissima.}
Thus $\erz X$ and $\zre X$ are both open in $\mc O_\down(K_n,\le)$
and hence elements of $\Fco{n}$. (Note that $X_\down$ and $X^\up$ are 
the closures of $X$ in the topologies $\mc O^\up(K_n,\le)$ and 
$\mc O_\down(K_n,\le)$ respectively.) 
We call a set $\Erz w$ for $w \in K_n$ a \emph{principal set} and a
set $\Zre w$ a \emph{co-principal set} (see Figure~\ref{f2}).

\begin{figure}[ht] 
\begin{picture}(6,4)
\put(1,0){\line(-1,5){0.75}}
\put(5,0){\line(1,5){0.75}}
\put(1,0){\line(1,0){4}}
\multiput(2.6,0.48)(0.15,0){11}{.}
\multiput(2.85,0.98)(0.15,0){8}{.}
\multiput(3.1,1.48)(0.15,0){5}{.}
\put(0.9,0.5){\line(1,0){1.65}}
\put(0.8,1){\line(1,0){2}}
\put(0.7,1.5){\line(1,0){2.35}}
\put(4.3,0.5){\line(1,0){0.8}}
\put(4,1){\line(1,0){1.2}}
\put(3.75,1.5){\line(1,0){1.55}}
\put(0.6,2){\line(1,0){4.8}}
\put(0.5,2.5){\line(1,0){5}}
\put(3,3){$\vdots$}
\put(3.3,1.55){$w$}
\put(3.4,1.8){\circle*{0.1}}
\multiput(2.325,0)(0.05,0.1){20}{-}
\multiput(4.375,0)(-0.05,0.1){20}{-}
\multiput(2.5,0)(0.1,0){19}{\line(0,1){0.1}}
\thicklines
\put(2.28,0.01){\line(1,2){1}}
\put(2.3,0){\line(1,2){1}}
\put(3.28,2.02){\line(1,0){0.24}}
\put(3.5,2){\line(1,-2){1}}
\put(3.515,2.01){\line(1,-2){1}}
\put(2.3,0.02){\line(1,0){2.2}}
\put(2.95,0.65){$\Erz w$}
\end{picture} 
\hfill
\begin{picture}(6,4)
\put(1,0){\line(-1,5){0.75}}
\put(5,0){\line(1,5){0.75}}
\put(1,0){\line(1,0){4}}
\multiput(0.95,0.48)(0.15,0){28}{.}
\multiput(0.85,0.98)(0.15,0){29}{.}
\multiput(0.75,1.48)(0.15,0){14}{.}
\multiput(0.65,1.98)(0.15,0){13}{.}
\multiput(0.55,2.48)(0.15,0){12}{.}
\multiput(3.6,1.48)(0.15,0){12}{.}
\multiput(3.8,1.98)(0.15,0){11}{.}
\multiput(4,2.48)(0.15,0){10}{.}
\put(3.03,1){\line(1,0){0.24}}
\put(2.8,1.5){\line(1,0){0.7}}
\put(2.55,2){\line(1,0){1.2}}
\put(2.3,2.5){\line(1,0){1.7}}
\put(3.15,3){$\vdots$}
\put(3.08,1.3){$w$}
\put(3.15,1.2){\circle*{0.1}}
\multiput(1.15,3)(-0.03,0.15){3}{.}
\multiput(5,3)(0.06,0.15){3}{.}
\multiput(0.955,0)(-0.02,0.1){36}{-}
\multiput(4.945,0)(0.02,0.1){36}{-}
\multiput(2.9,1)(-0.05,0.1){26}{-}
\multiput(3.3,1)(0.05,0.1){26}{-}
\multiput(1.1,0)(0.1,0){39}{\line(0,1){0.1}}
\multiput(3,0.9)(0.1,0){4}{\line(0,1){0.1}}
\thicklines
\put(0.98,0){\line(-1,5){0.75}}
\put(5.02,0){\line(1,5){0.75}}
\put(1,0.02){\line(1,0){4}}
\put(3.03,1){\line(-1,2){1.35}}
\put(3.05,1){\line(-1,2){1.35}}
\put(3,1.02){\line(1,0){0.26}}
\put(3.02,1){\line(1,0){0.25}}
\put(3.25,1){\line(1,2){1.35}}
\put(3.27,1){\line(1,2){1.35}}
\put(1,1.62){$\Zre w$}
\end{picture} 
\caption{A principle set $\Erz w$ and a co-principle set $\Zre w$\label{f2}} 
\end{figure}

The following theorem and its corollary are main results of Bellissima%
\footnote{Lemma~2.6, Theorem~2.7 and Corollary~2.8 in \cite{B}; our
Corollary~\ref{Cor-Thm-B2}\,(b) is implicit in Bellissima's Lemma~2.6.}. 
We follow his proof except that we simplify notations and arguments
and that we get shorter formulae.%
\footnote{This is mainly because Bellissima's $\phi_1$ is implied by
  the second conjunct of his $\phi_2$.} 
Of course, an empty disjunction stands for the formula $\bot$ and an
empty conjunction for $\top$.

\begin{Theorem} \label{Thm-B2} 
  For every $w \in K_n$, there are formulae $\psi_w$ and $\psi'_w$ such
  that $[\![\psi_w]\!] = \Erz{w}$ and $[\![\psi'_w]\!] = \Zre{w}$. They
  can be defined by induction on the foundation rank of $w$ as follows. 
  If\/ $Y_{\max}$ denotes the maximal elements of $Y$, then 
\begin{align*}
  \psi_{w_{\beta,Y}} &:=\ \Big(\Big( \bigvee_{w \rin Y_{\max}} \!\!\psi'_{w}
  \;\lor \bigvee_{P_i \notin \beta} P_i \Big) \;\to\! \bigvee_{w \rin Y_{\max}}
  \!\!\psi_w \Big) \;\land\; \bigwedge_{P_i \in \beta} P_i \\ 
  \psi'_{w_{\beta,Y}} &:=\quad \psi_{w_{\beta,Y}} \;\to\! \bigvee_{w \rin Y_{\max}} \!\!\psi_{w} 
\end{align*}
It follows that if $w$ has foundation rank $d$, then the implication
depth of $\psi_w$ is at most $2d+1$ and of $\psi'_w$ at most $2d+2$.
\end{Theorem} 

\PROOF 
Let $w_{\beta,Y}$ be of foundation rank $d+1$. We assume the
proposition to be shown for all points of smaller foundation rank in
$K_n$, and conclude by induction. 

By induction, $[\![\bigvee_{w \rin Y_{\max}} \psi_{w}]\!] = Y$ 
(if $d=-1$, then $Y = \emptyset$, and everything works as well). 
Now $\psi_{w_{\beta,Y}}$ is intuitionistically equivalent to a
conjunction of three formulae that define the following subsets 
of $K_n$: 

\begin{align*} 
A_1 \ :=\ 
[\![ \bigvee_{w \rin Y_{\max}} \!\!\psi'_{w} \;\to\! \bigvee_{w \rin Y_{\max}} \!\!\psi_w]\!]
\ &=\ \Big\{ w \Bigm| \forall v \le w \; 
      \Big(\, v \in \!\! \bigcap_{z \rin Y_{\max}} \!\! [\![ \psi'_{z} ]\!]\cpl
      \ \cup\!\! \bigcup_{z \rin Y_{\max}} \!\! [\![ \psi_z]\!] \,\Big) \,\Big\} \\
\ &=\ \big\{ w \bigm| \forall v \le w \; \big(\, Y \subseteq \Erz{v} 
      \ \text{ or }\ v \in Y \,\big) \,\big\} \\[1ex]
\ &\subseteq\ \phantom{\big\{} B_1 \ := \ 
      \big\{ w \bigm| \Erz w \cap K_n^d = Y \big\} \;\cup\; Y \\[2ex]
A_2 \ :=\ 
[\![ \bigvee_{P_i \notin \beta} P_i \;\to\! \bigvee_{w \rin Y_{\max}} \!\!\psi_w]\!]
\ &=\ \Big\{ w \Bigm| \forall v \le w \; 
      \Big(\, v \in \! \bigcap_{P_i \notin \beta} \! [\![ P_i ]\!]\cpl
      \ \cup\!\! \bigcup_{z \rin Y_{\max}} \!\! [\![ \psi_z]\!] \,\Big) \,\Big\} \\
\ &=\ \big\{ w \bigm| \forall v \le w \; \big(\, \val(v) \subseteq \beta 
      \ \text{ or }\ v \in Y \,\big) \, \big\} \\[2ex]
\text{ and }\qquad 
A_3 \ :=\ 
[\![ \bigwedge_{P_i \in \beta} P_i ]\!] \ &=\ 
  \big\{w \in K_n \bigm| \beta \subseteq \val(w) \big\}. 
\end{align*} 
One sees that $Y\subseteq A_i \cap K_n^d$ for all $i$, and
$A_1 \cap K_n^d = Y$. Thus $[\![\psi_{w_{\beta,Y}}]\!] \cap K_n^d = Y$. 
Moreover, if $w \in [\![\psi_{w_{\beta,Y}}]\!] \setminus K_n^d$, then
$w$ has the following property: for all $v \le w$, either $v \in Y$, or
$\val(v) = \beta$ (from $A_2$ and $A_3$) and $\Erz v \cap K_n^d = Y$
(from $B_1$). By construction of $\mf K_n$, there is only one such
point, namely $w_{\beta,Y}$.

Then $[\![\psi'_{w_{\beta,Y}}]\!]$ is by definition the largest
decreasing set contained in $[\![\psi_{w_{\beta,Y}}]\!]\cpl \cup 
[\![\bigvee_{w \rin Y_{\max}} \psi_{w}]\!] = K_n \setminus\{w_{\beta,Y}\}$, 
which is exactly $\Zre{w_{\beta,Y}}$. 
\END

\begin{Cor} \label{Cor-Thm-B2} 
\begin{enumerate}[{\bf(a)}]
\item[]
\item 
  The principal and the co-principle sets are in $\Fr{n}$, hence also
  all finite sets. 
\item 
  For $w \in K_n$, if $a = \Erz w$ and $b = \Zre w$, then $b = a \hto a^-$.
\end{enumerate} 
\end{Cor} 

\PROOF
(b) follows because ${\Erz{w_{\beta,Y}}}^- = Y = 
 [\![\bigvee_{w \rin Y_{\max}} \psi_{w}]\!] $. 
\END

\medskip
From the construction in Theorem~\ref{Thm-B1} we will mainly use two
properties: 
The ``filtration'' of the Kripke model into levels of finite foundation 
rank. And the property that any finite set of at least two incomparable 
elements in $K_n$ has a common successor without other predecessors. 
Moreover, we need the following result from \cite{B}:

\begin{Fact}[Theorem 3.0 in \cite{B}] \label{F-Bneu} 
\begin{enumerate}[{\bf(a)}]
\item[]
\item 
  The principal sets are exactly the $\Hcup$-irreducible elements of
  both algebras, $\Fr{n}$ and $\Fco{n}$. 
\item 
  The co-principal sets are exactly the $\hcap$-irreducible
  sets of both algebras. 
\end{enumerate} 
\end{Fact} 

The theorem in \cite{B} is formulated for $\Fr{n}$ only, but the
proof works as well for $\Fco{n}$. For the sake of completeness, we
add a sketch of the proof: 

\PROOF 
(a) For $X \in \Fco{n}$ we have $X = \bigcup_{w \in X} \Erz{w}$. It
follows that $X$ is $\Hcup$-irreducible, if and only if there is a
greatest element $w_0$ in $X$, if and only if $X = \Erz{w_0}$. 
 
(b) If there are two minimal elements $w_0,w_1 \in K_n \setminus X$,
then $X = (X \cup \{w_0\}) \cap (X \cup \{w_1\})$ is not
$\hcap$-irreducible. Conversely, $\Zre w$ has a unique successor,
namely $\Zre w \cup \{w\}$, hence $\Zre w$ is $\Hcap$-irreducible. 
\END

\medskip 
Let $\Fkl{n}$ be the sub-Heyting algebra of $\Fr{n}$ generated by all 
$\Hcup$-irreducible elements of $\Fr{n}$. Thus $\Fkl{n} = \mathbf{B}_n^{}$ 
in Bellissima's notation in \cite{B}.

\begin{Fact}[Theorem~4.4 in \cite{B}] \label{F-Fsmile} 
  For $n>1$, the algebra $\Fkl{n}$ is not finitely generated, because
  the sub-algebra generated by all $\Erz{w}$ for $w$ of foundation
  rank $\le d$ can't separate points of $K_n$ of higher foundation
  rank that differ only by their valuations.
\end{Fact}

In particular, it follows that $\Fkl{n}$ is not isomorphic to $\Fr{n}$
for $n>1$. To our knowledge, Grigolia has shown that no proper
sub-algebra of $\Fr{n}$ is isomorphic to $\Fr{n}$. 
Without being explicitly mentioned, it is clear in \cite{B} that 
$\Fkl{1} = \Fr{1} = \Fco{1}$ and that $\Fr{n} \neq \Fco{n}$ for
$n>1$.

\begin{Fact}[Lemma~4.1 in \cite{B}] \label{F-anti} 
  For $n>1$, there is an infinite antichain in $K_n$.
\end{Fact}

\section{Some consequences} \label{Sec-csq} 

This section collects some rather immediate consequences from
Bellissima's construction that are not explicitly mentioned in
Bellissima's paper, and which we will use in our analysis of
Bellissima's setting. Other consequences are collected in 
section~\ref{Sec-letzt}. Many of the results hold in a much wider
context, see \cite{DJ2}.

\subsection{More on irreducible elements}

For $a$ in a Heyting algebra, let us define the \emph{supports} 
\begin{align*}
  \supp_{\Hcup}(a) \ :=\ & \big\{ x \text{ $\textstyle\Hcup$-irreducible} \bigm| x\hle a \big\}, \\
  \supp_{\Hcap}(a) \ :=\ & \big\{ x \text{ $\textstyle\Hcap$-irreducible} \bigm| a\hle x \big\}. \\
  \suppmin_{\Hcap}(a) \ :=\ & \;\text{the minimal elements in } \supp_{\Hcap}(a) 
\end{align*}

\begin{Lemma} \label{L-irr}
\begin{enumerate}[{\bf(a)}]
\item[]
\item 
  The $\hcap$-irreducible elements are $\Hcap$-irreducible in both,
  $\Fr{n}$ and $\Fco{n}$. 
\item 
  For each $a\in \Fco{n}$, we have 
  $$a\ =\ \bighcup\; \supp_{\Hcup}(a) 
     \ =\ \bighcap\; \suppmin_{\Hcap}(a),$$
  and $a \neq \bighcap S$ for every proper subset $S \subset
  \suppmin_{\Hcap}(a)$. 
\end{enumerate} 
\end{Lemma}

It follows from part (a) of this Lemma and from Fact~\ref{F-Bneu}\,(a) 
that the (well-founded) partial ordering of $K_n$ equals the partial
ordering of the $\Hcap$-irreducibles, and that of the $\Hcup$-irreducibles. 
In particular, the (co-)foundation rank of $w$ in $K_n$ equals the
\mbox{(co-)}foundation rank of $\Erz{w}$ in the partial ordering of the
$\Hcup$-irreducibles and also the (co-)foundation rank of $\Zre w$ in
the partial ordering of the $\Hcap$-irreducibles. Hence all these
foundation ranks are finite, and all these co-foundation ranks equal
to $\infty$, due to the existence of an infinite chain 
$w < w' < w'' < \cdots$ in $K_n$. Foundation and co-foundation ranks
are used to compute dimensions and co-dimensions, cf. \cite{DJ2}
remark~6.9.

\PROOF 
(a) has already be shown in the proof of Fact~\ref{F-Bneu}\,(b). 

(b) For any decreasing set $a \subseteq K_n$, we have $a = 
\bigcup_{w \in a} \Erz w = \bigcap_{w \notin a} \Zre w$, which proves
the first equality and the second for the full $\supp_{\Hcap}$. But
because the order is well-founded on the $\Hcap$-irreducible elements,
it is enough to keep the minimal elements in the support. The last
statement is clear by definition. 
\END

\begin{Cor} \label{Cor-F-Be} 
\begin{enumerate}[{\bf(a)}]
\item[]
\item 
  The $\Hcup$-irreducibles and $\Hcap$-irreducibles of $\Fco{n}$ and
  of $\Fr{n}$ are the same. 
\item
  Any element of $\Fco{n}$ is a (possibly infinite) union as well as a
  (possibly infinite) intersection of elements of $\Fr{n}$.
\end{enumerate}
\end{Cor}

A characterisation of the $\hcup$-irreducible elements can be found in 
Section~\ref{S-hcupirr}.

\begin{Ex} 
Lemma~\ref{L-irr} doesn't work with the maximal elements in
$\supp_{\Hcup}$, because the order on the $\Hcup$-irreducible elements
is not anti-well-founded. For example, if $n \ge 2$, then there are no
maximal elements in $\supp_{\Hcup} \big( [\![P_i]\!] \big)$. 
However, for $\Fco{n}$ Corollary~6.10 in \cite{DJ2} provides a
decomposition into $\hcup$-irreducible components. 

\PROOF 
Say $i=2$. One can show that for each $k$, there are at least two
elements in $K_n$ of foundation rank $k$ and with $P_2$ in the
valuation. For example because there is a copy of the Kripke model
$\mf K_1$ inside the points $w$ of $K_n$ with $P_2 \in \val(w)$ 
---just add $P_2$ to the valuations of the points of $\mf K_1$,
i.e. start with the points $w_{\{P_2\},\0}$ and $w_{\{P_1,P_2\},\0}$--- 
and $\mf K_1$ is well known to have two points of each foundation rank.
Now if $w \in [\![P_2]\!]$, choose $v \in [\![P_2]\!]$, $v\neq w$, of
same foundation rank. Then the point $w_{\{P_2\},\Erz{v,w}}$ shows
that $\Erz w$ is not maximal in $\supp_{\Hcup} \big( [\![P_2]\!] \big)$.
\END
\end{Ex}

\begin{Rem} \label{Rem-suppKripke}
We can identify the underlying partially ordered set of the Kripke
model $\mf K_n$ with either the $\Hcup$-irreducibles via $w \mapsto
\Erz{w}$, or with the $\Hcap$-irreducibles via $w \mapsto \Zre{w}$,
both with the order induced by the partial order $\hle$ of the Heyting
algebra $\Fr{n}$ (which is one of the reasons to work with the
reversed order on the Kripke model). 
We can further identify the $\Hcap$-irreducibles with the principal
prime ideals of $\Fr{n}$ that they generate, ordered by inclusion. 
By Lemma~\ref{L-irr}\,(a), all principal prime ideals are of that
form. The valuations of the Kripke model can also be recovered from
$\Fr{n}$, as will be explained in Section~\ref{Sec-vier}. 

A survey about the relationship between Kripke models and Heyting algebras
can be found in the doctoral thesis of Nick Bezhanishvili \cite{Nick}. 
\end{Rem}

\begin{Rem} \label{Rem-suppminneu}
The decomposition $a = \bighcap\; \suppmin_{\Hcap}(a)$ of
Lemma~\ref{L-irr} provides a sort of infinite (conjunctive) normal
form for elements of $\Fco{n}$. Provided the partial order on the
$\Hcap$-irreducibles is known, the Heyting algebra operations can be
computed from the supports in a first order way: First note that
$\suppmin_{\Hcap}$ and $\supp_{\Hcap}$ can be computed from each
other, and then
\begin{align*}
\suppmin_{\Hcap}(a\hcup b) &\ =\ \text{the minimal elements in }
  \;\supp_{\Hcap}(a) \,\cap\, \supp_{\Hcap}(b) \\\
\suppmin_{\Hcap}(a\hcap b) &\ =\ \text{the minimal elements in }
  \;\supp_{\Hcap}(a) \,\cup\, \supp_{\Hcap}(b) \\
\suppmin_{\Hcap}(a\hto b) &\ =\ \text{the minimal elements in }
  \;\supp_{\Hcap}(b) \,\setminus\, \supp_{\Hcap}(a) 
\end{align*}
Note that every set of pairwise incomparable $\Hcap$-irreducible
elements forms the $\suppmin_{\Hcap}$ of an element of $\Fco{n}$,
namely of its intersection. 
\end{Rem}

\begin{Qu}
  Can we characterise those sets that correspond to elements of $\Fr{n}$?
\end{Qu}

Because $\Erz w \in \supp_{\Hcup}(a) \iff \Zre w \notin
\supp_{\Hcap}(a)$, one can translate the above rules into computations
of the Heyting algebra operations from the $\Hcup$-supports, but they
are less nice.

\subsection{The completion as a profinite limit}

Consider in $\Fr{n}$ for each $i$ the filter $(K_n^i)^\up =
\{a\in\Fr{n} \mid K_n^i \subseteq a\}$ generated by $K_n^i$, and
denote the corresponding congruence relation by $\equiv_i$. Thus
$$a \equiv_i b \ \iff \ a \cap K_n^i = b \cap K_n^i.$$ 
We denote the quotient $\Fr{n}/{\equiv_i}$ by $\Fr[i]{n}$ and
the canonical epimorphism by $\pi_i$. It extends to an epimorphism
$\wh{\pi}_i: \Fco{n} \to \Fr[i]{n}$, where $\wh{\pi}_i^{-1}(1) = 
\{ a \in \Fco{n} \mid K_n^i \subseteq a\}$ is the filter generated by
$K_n^i$ in $\Fco{n}$. For simplicity, we denote it and the
corresponding congruence relation again by $(K_n^i)^\up$ and
$\equiv_i$ and we will identify $\Fco{n}/{\equiv_i}$ with $\Fr[i]{n}$. 

$\Fr[i]{n}$ is naturally isomorphic to the finite Heyting algebra 
$\mc O_\down(K_n^i,\le)$ via ``truncation'' $\pi_i(a) \mapsto a \cap
K_n^i$ (the surjectivity needs part (a) of Corollary~\ref{Cor-F-Be}).
Again, we will identify both without further mentioning.%
\footnote{There is a certain ambiguity here that should not harm: 
  As $K_n^i$ is also an element of $\Fr{n}$, there is a map $\Fr{n}
  \to \Fr{n}$, $a \mapsto a \cap K_n^i$. The image of this map is a
  sub-poset (but not a sub-algebra) of $\Fr{n}$ isomorphic to $\Fr[i]{n}$.}

\begin{Rem} 
There is a natural notion of dimension (more precisely: dual
codimension, see \cite{DJ2} and \cite{D-preprint}) in lattices such 
that $\Fr[i]{n}$ is the free Heyting algebra over $n$ generators of
dimension $i$. Among universal algebraists, it is known as the free
algebra generated by $n$ elements in the variety of Heyting algebras
satisfying $P_{i+1} = 1$ (where $P_{i+1}$ is defined inductively as 
$x_{i+1} \hcup (x_{i+1} \to P_i)$ with $P_0 = 0$). 
\end{Rem}

\begin{Prop} \label{P-rf}
  $\Fr{n}$ and $\Fco{n}$ are residually finite, \textit{i.e.} for each
  element $a \neq 1$ there is a homomorphism $\psi$ onto a finite
  Heyting algebra such that $\psi(a) \neq 1$. Moreover, $\psi$ can be
  chosen to be some $\pi_i$ or $\wh{\pi_i}$ respectively.
\end{Prop}

\PROOF 
This is clear from the above considerations and the construction of
$K_n$ as union of the $K_n^i$. 
\END

Because $K_n^i \subseteq K_n^j$ for $i < j$, the morphism $\pi_i:
\Fr{n} \to \Fr[i]{n}$ factors through $\Fr[j]{n}$ for $i < j$ and
yields an epimorphism $\pi_{ji}: \Fr[j]{n} \to \Fr[i]{n}$. The system
of maps $\pi_{ji}$ is compatible, hence the projective limit
$\varprojlim\Fr[i]{n}$ exists. 
Because of the universal property of the projective limit and because
the system of morphisms $\wh{\pi_i}: \Fco{n} \to \Fr[i]{n}$ is
compatible as well, there is a natural morphism $\Fco{n} \to
\varprojlim\ \Fr[i]{n}$. 

$$\begin{array}{ccccr@{}c@{}c@{}c@{}clcc}
  \Fco{n} && \hookleftarrow &&&& \Fr{n} \\[1ex] 
  \downarrow &&& \cdots & {}^{\pi_{i+1}} & \sswarrow 
  & \multicolumn{1}{r@{\,}}{{}^{\pi_i}} & \ssearrow && \dotsm \\[1ex] 
  \varprojlim\,\Fr[i]{n} & \longonto & \dotsm & \longonto &\Fr[i+1]{n} 
  && \underset{\pi_{i+1,i}}{\loongonto} && \Fr[i]{n} & \longonto & \dotsm 
\end{array}$$

\begin{Prop} 
$\Fco{n}$ is the projective limit of the finite Heyting algebras
  $\Fr[i]{n}$, in symbols 
$$\Fco{n} \ = \ \underset{i\in\mbb N}{\varprojlim}\ \Fr[i]{n}.$$ 
Any finite quotient of $\Fr{n}$ factors through some $\Fr[i]{n}$,
hence $\Fco{n}$ is the profinite completion of $\Fr{n}$, i.e. the
projective limit of all finite epimorphic images of $\Fr{n}$. 
\end{Prop} 

The content of this proposition or parts of it is, in various forms,
contained in several articles. To our knowledge, Grigolia was the
first to embed $\Fr{n}$ into ${\varprojlim}\,\Fr[i]{n}$, see
\cite{Grigolia-buch} or \cite{Grigolia-Aut}. 
More about profinite Heyting algebras can be found in \cite{BGMM} and
\cite{Bz2}. In the first of these two articles, the profinite
completion of a Heyting algebra is embedded in the topological Heyting
algebra of its dual space, i.e. its prime spectrum. 

\PROOF 
The natural morphism $\Fco{n} \to \varprojlim\ \Fr[i]{n}$ is in fact
the map $a \mapsto (\wh{\pi_i}(a))_{i \in \mbb N} = 
(a \cap K_n^i)_{i \in \mbb N}$. It is injective by
Proposition~\ref{P-rf} and it is surjective because for any compatible
system $(a_i)_{i \in \mbb N}$ of decreasing sets $a_i \subseteq
K_n^i$, the union $\bigcup_{i \in \mbb N} a_i$ is a pre-image in 
$\Fco{n}$. 

Now let $\Phi$ be a filter such that $\Fr{n}/{\equiv_\Phi}$ is finite. 
Choose an $i$ such that for all elements in $\Fr{n}/{\equiv_\Phi}$ there 
is a pre-image in $\Fr{n}$ with pairwise distinct images in $\Fr[i]{n}$. 
Then $\Fr{n} \to \Fr{n}/{\equiv_\Phi}$ factors through $\Fr[i]{n}$. 
\END

The epimorphism $\pi_{i+1,i}: \mc O_\down(K_n^{i+1},\le) \to \mc O_\down(K_n^i,\le)$
is induced from the inclusion $(K_n^i,\le) \into (K_n^{i+1},\le)$. The
theorem shows that $\mc O_\down$ behaves like a contravariant functor, i.e. 
$$\Fco{n} \;=\; \mc O_\down(K_n,\le) \;=\; 
  \mc O_\down \big(\underset{i\in\mbb N}{\varinjlim}\ (K_n^i,\le)\big) \;=\;  
  \underset{i\in\mbb N}{\varprojlim}\ \mc O_\down (K_n^i,\le) \;=\; 
  \underset{i\in\mbb N}{\varprojlim}\ \Fr[i]{n}$$

Note that any profinite lattice is a complete lattice as projective
limit of complete lattices, and that any profinite structure is well
known to be a compact Hausdorff topological space (the topology being
the initial topology for the projections onto the finite quotients in
the defining system, equipped with the discrete topology). We are
going to examine this topology a little further:

\medskip\textbf{Definition } 
For $a,b \in \Fco{n}$, define their distance to be 
$$\mr d(a,b) 
 \;:=\; 2^{- \min\{i \,\mid\, \wh{\pi}_i(a) \neq \wh{\pi}_i(b) \} } 
 \; =\; 2^{- \min\{i \,\mid\, a \cap K_n^i \,\neq\, b \cap K_n^i \} }$$ 
with the convention $2^{-\min\emptyset} = 0$.

\begin{Theorem} \label{Thm-metric} 
  $(\Fco{n}, \mr d)$ is a compact Hausdorff metric space. The metric
  topology is the profinite topology, and the Heyting algebra
  operations are continuous. $\Fco{n}$ is the metric completion of its
  dense subset $\Fr{n}$. 
\end{Theorem} 

This is analogous to properties of the $p$-adic integers $\mbb Z_p$ as 
a projective limit of the rings $\mbb Z/p^n\mbb Z$. See also
Theorem~6.1 in \cite{DJ2} for a generalisation to a wider class of
Heyting algebras comprising the finitely presented ones.

\PROOF 
It is straightforward to check that $\mr d$ is a metric: 
symmetry and the ultrametric triangular inequality $\mr d(a,c) \le
\max\big\{\mr d(a,b),\mr d(b,c)\big\}$ are immediate from the
definition, and $\mr d(a,b) = 0 \Longleftrightarrow a=b$ holds by 
Proposition~\ref{P-rf}. 

The profinite topology and the metric topology coincide because they
have the same basis of (cl-)open sets $\wh{\pi_i}^{-1}(a) = 
\{ x \mid \mr d(x,x_0) < 2^{-i+1}\}$  for $a \in \Fr[i]{n}$ and $x_0$
with $\wh{\pi_i}(x_0)=a$. 

It is a standard result that a projective limit of compact Hausdorff 
spaces is again compact and Hausdorff. (Here, this is easily seen
directly: Metric topologies are always Hausdorff. If a family of
closed sets with 
the finite intersection property is given, we may suppose the closed
sets to be of the form $\wh{\pi_i}^{-1}(a_i)$. The fip implies that
every $i$ appears only once and that the $a_i$ are compatible. Hence
$(a_i)_{i\in\omega}$ is an element of $\varprojlim\Fr[i]{n}$ in the
intersection of the family and the topology is compact.)

By definition of the profinite topology, the maps $\wh{\pi_i}$ are 
continuous. Then the continuity of the Heyting algebra operations on
the $\Fr[i]{n}$ (with discrete topology!) lifts to $\Fco{n}$, e.g.\
$\hcap^{-1}\big(\wh{\pi_i}^{-1}(a)\big) = 
(\wh{\pi_i}\times\wh{\pi_i})^{-1}\big(\hcap^{-1}(a)\big)$ is open. 
An element $a\in\Fco{n}$ is a limit of the sequence 
$(a\cap K_n^i)_{i\in\omega}$ in $\Fr{n}$. Therefore $\Fco{n}$ is the
metric completion and thus $\Fr{n}$ is dense in $\Fco{n}$. 
\END

As points are closed, we get that any term in the language of Heyting
algebras (even with parameters) defines a continuous map.

\begin{Rem}
There is a fundamental difference between the cases $n=1$ and $n>1$. 
In case $n=1$, the map $\pi_{i+1\,i}: \Fr[i+1]{1} \to \Fr[i]{1}$ has a
kernel of $3$ elements, but is otherwise injective. It follows that
the map $\pi_i: \Fco{1} \to \Fr[i]{1}$ is injective on $\Fco{1}
\setminus (K_n^i)^\up$, and that $\Fr{1} = \Fco{1}$. In case $n>1$,
the size of $\Fr[i]{n}$ is growing faster than exponentially with $i$.
Moreover, the maps $\pi_{i+1\,i}: \Fr[i+1]{n} \to \Fr[i]{n}$ are non-injective
enough to allow a tree of elements $a_{s_1,\dots,s_i} \in \Fr[i]{n}$
with $s_j \in \{0,1\}$ such that $\pi_{i+1\,i}(a_{s_1,\dots,s_{i+1}})
= a_{s_1,\dots,s_i}$. Hence $\Fco{n}$ has size continuum and differs
from $\Fr{n}$. 
\end{Rem}

\subsection{Dense sub-algebras}

\textbf{Definition }
We say that $H$ is \emph{dense in $\Fco{n}$} if $H$ is a dense Heyting
sub-algebra of $\Fco{n}$ in the metric topology. By Theorem~\ref{Thm-metric}, 
the metric topology equals the profinite topology, hence density of a
sub-algebra $H$ of $\Fco{n}$ means that all the induced maps $\pi_d:H
\to \Fr[d]{n}$ are surjective.

\begin{Lemma} \label{L-isolated} 
  The isolated points in $\Fco{n}$ are exactly the finite elements. 
  They are dense in $\Fco{n}$. 
\end{Lemma} 

\PROOF
If $a$ is finite $\subseteq K_n^i$, then there is no other element in
the $2^{-(i+2)}$-ball around $a$, hence $a$ is isolated. If $a$ is
infinite, then $a$ is the limit of a sequence of finite elements
distinct from $a$, namely $a = \lim_{i\in\omega} (a\cap K_n^i)$. 
Hence $a$ is not isolated, but a limit of isolated points.
\END

\medskip
Recall that $\Fkl{n}$ is the sub-Heyting algebra of $\Fr{n}$ generated
by all $\Hcup$-irreducible elements of $\Fr{n}$, i.e. the sub-algebra
generated by all finite elements. Thus Lemma~\ref{L-isolated}
immediately implies:

\begin{Prop} \label{Pr-Fkl}
  $\Fkl{n}$ is the smallest dense sub-Heyting algebra of $\Fr{n}$.
\end{Prop}

As we have mentioned in Fact~\ref{F-Fsmile}, Bellissima has shown that
$\Fkl{n}$ is not isomorphic to $\Fr{n}$ for $n\geqslant 2$. If one
knew that $\Fkl{n} \neq \Fr{n}$, then Proposition~\ref{Pr-Fkl} would
yield an easier proof, because $\Fr{n}$ then has a proper dense 
sub-algebra, namely $\Fkl{n}$, whereas $\Fkl{n}$ does not. We will see
later that $\Fkl{n} \not\equiv \Fr{n}$ for $n\geqslant 2$.

\begin{Lemma} \label{L-denseirr} 
  If $H$ is dense in $\Fco{n}$, then $H$ has the same
  $\Hcup$-irreducibles and the same $\Hcap$-irreducibles as $\Fco{n}$. 
\end{Lemma} 

\PROOF 
Because of the density, all the finite sets are in $H$. Thus in
particular all the principal sets are in $H$, and as they are
$\Hcup$-irreducible in $\Fco{n}$, they remain irreducible in $H$. 
If $X \in H$ is not principal, then $X$ is the proper union of its
$\supp_{\Hcup}$, hence not $\Hcup$-irreducible.

Because of Corollary~\ref{Cor-F-Be}\,(b), $H$ also contains all the
$\Hcap$-irreducibles of $\Fco{n}$, and they remain for trivial reasons
$\Hcap$-irreducible in $H$. With the same argument as above, one sees
that there are no other $\Hcap$-irreducible elements in $H$, because
any other element if the proper intersection of its $\supp_{\Hcap}$. 
\END

\medskip 
It follows that $\Fkl{n}$ is also the sub-Heyting algebra generated 
by the $\Hcap$-irreducible elements of $\Fco{n}$, because any finite
set is the intersection of finitely many $\Hcap$-irreducible elements,
namely its $\suppmin_{\Hcap}$. See also section~\ref{S-approx} for
more on $\Fkl{n}$, and Proposition~6.12 of \cite{DJ2} for a
generalisation of $\Fkl{n}$ to, among others, finitely presented
Heyting algebras.

\begin{Rem} 
If $n>1$ and $b_1, b_2$ are two incomparable $\Hcap$-irreducible
elements in $\Fr{n}$, then $\suppmin_{\Hcap}(b_1 \hcup b_2)$ is
infinite. Thus not all elements of $\Fkl{n}$ have finite
$\suppmin_{\Hcap}$. This is a reason why an explicit description of
$\Fkl{n}$ is not as easy as one could first think.
\end{Rem} 

\PROOF
If $b_i = \Zre{w_i}$, then $\Zre{w_{\emptyset,\Erz{w_1,w_2,v}}}$ is in
$\suppmin_{\Hcap}(b_1 \hcup b_2)$ for any $v$ that is incomparable
with $w_1$ or with $w_2$. If $n>1$, there are infinitely many such $v$. 
\END

\begin{Rem} \label{Rem-section} 
There are two canonical sections of the projection map $\pi_i: \Fco{n}
\to \Fr[i]{n}$: 
\begin{align*} 
  \text{the \emph{minimal section}}\quad \sigma_i^{\min}:&  
  \ x \;\mapsto\; \bighcap\; \big\{ y \in \Fco{n} \bigm| \pi_i(y) = x \big\} \\ 
  \text{ and the \emph{maximal section}}\quad \sigma_i^{\max}:&  
  \ x \;\mapsto\; \bighcup\; \big\{ y \in \Fco{n} \bigm| \pi_i(y) = x \big\}
\end{align*}
In fact, both have images in $\Fkl{n}$: each $\sigma_i^{\min}(x)$ is a
finite set, thus a finite union of $\Hcup$-irreducibles, and each
$\sigma_i^{\max}(x)$ is what one might call a ``co-finite set'',
namely a finite intersection of $\Hcap$-irreducibles.  If one sees
$x\in \Fr[i]{n}$ as a subset of the Kripke model $\mf K^i_n$ for
$\Fr[i]{n}$, then the minimal section maps $x$ to itself but now seen
as a subset of the Kripke model $\mf K_n$. 
One can compute $\sigma_i^{\min}(x) = K_n^i \hcap \sigma_i^{\max}(x)$ 
and $\sigma_i^{\max}(x) = (K_n^i \hto \sigma_i^{\min}(x))$, and one can
check that $\sigma_i^{\min}$ is a $\{0,\hcap,\hcup,\hle\}$-homomorphism 
and $\sigma_i^{\max}$ is a $\{0,1,\hcap,\hto,\hle\}$-homomorphism. 

The quotient $\Fr[0]{n}$ of $\Fr{n}$ is isomorphic to $\mf P(K_n^0)$,
which, via the valuation, can be identified with the free Boolean
algebra $\mf P(\mf P(P_1,\dots,P_n))$ over the free generators 
$P_1,\dots,P_n$. The image of the maximal section of $\pi_0$ consists
exactly of the regular elements of $\Fr{n}$. Thus the regular elements
form (as a sub-poset, but not as a sub-algebra) a free Boolean algebra
over $n$ generators. (This is easy to see with Bellissima's characterisation
of regular elements in Corollary 2.8 of \cite{B}, and much of it is
already in \cite{McTarski}.) 
\end{Rem}

\section{Reconstructing the Kripke model} \label{Sec-vier}

\subsection{Another duality}

The following might be well known in lattice theory. 
In a (sufficiently) complete lattice, define 
$$a^\hcap\ :=\ \bighcup\, \{ x \mid a \not\hle x \} \quad\text{ and }\quad
a^\hcup\ :=\ \bighcap\, \{ x \mid x \not\hle a \}$$

\begin{Lemma} \label{L-dual}
Suppose $\Lambda$ is a complete lattice that satisfies both 
infinite distributive laws. Then the maps $a \mapsto a^\hcap$ and 
$b \mapsto b^\hcup$ are inverse order-preserving bijections between
the $\Hcup$-irreducibles and the $\Hcap$-irreducibles. Moreover, 
$a \not\hle a^\hcap$ and $b^\hcup\not\hle b$.
\end{Lemma} 

\PROOF 
The maps are order-preserving by definition (and the transitivity of $\hle$).
Let $a$ be $\Hcup$-irreducible and suppose $a \hle a^\hcap$. Then 
$a = a \hcap a^\hcap = \bighcup \{ a \hcap x \mid a \not\hle x \}$. 
Then the $\Hcup$-irreducibility of $a$ implies $a = a \hcap x$ for some
$x \not\hge a$: contradiction. Hence $a^\hcap$ is the greatest element
$x$ with the property $a \not\hle x$. It follows that $a^\hcap \hleneq 
a \hcup a^\hcap$, and if $a^\hcap \hleneq b$, then $a \hle b$. Thus 
$a^\hcap \hcup a$ is the unique successor of $a^\hcap$, which therefore
is $\Hcap$-irreducible. Since dually $a^{\hcap\hcup}$ is the minimal
element $x$ with the property $x \not\hle a^\hcap$, and $a \not\hle a^\hcap$, 
we get $a^{\hcap\hcup} \hle a$. If we had $a \not\hle a^{\hcap\hcup}$,
then $a^{\hcap\hcup} \hle a^\hcap$ by definition of the latter. 
But dually to the argument above, $a^{\hcap\hcup}$ is the smallest
element $x$ with $x \not\hle a^\hcap$: contradiction and $a = a^{\hcap\hcup}$.
The remaining parts are by duality. 
\END

One can reformulate Lemma~\ref{L-dual} partially as 
$$\mc O_\down \Big(\text{$\Hcap$-irreducibles},\hle \Big) \ \isom\ 
  \mc O_\down \Big(\text{$\Hcup$-irreducibles},\hle \Big).$$

\begin{Rem} \label{R-Hcup->Hcap}
(1) As $x\notin\supp_{\Hcap}(a) \iff a \not\hle x \iff x^\hcup \hle a$ 
by definition of $x^\hcup$, it follows that 
$$\supp_{\Hcup}(a)\ =\ \big\{\, x^\hcup \bigm| x \text{ $\textstyle\Hcap$-irreducible}, 
  x \notin \supp_{\Hcap}(a) \big\}\ \text{ if $a \neq 0$},$$ 
and of course the dual statement also holds. 

(2) $\Fco{n}$ as a lattice of sets satisfies the hypotheses of 
Lemma~\ref{L-dual}. In that special situation, we have that any
$\Hcup$-irreducible $a$ is of the form $\Erz w$. The element $a^\hcap$
is then the corresponding co-principal set $\Zre w$, and we have seen 
in Corollary~\ref{Cor-Thm-B2}\,(b) that $a^\hcap = a \hto a^-$. 
For a $\Hcap$-irreducible $a$ however, we have $a^+ \hto a = a$. 
Considering $\Fco{n}$ with its co-Heyting structure, we get the
``dual'' rule $a^\hcup = a^+ - a$. (Recall from the Example~\ref{Ex}
that all partial orders are bi-Heyting.) 
\end{Rem}

\subsection{The spectrum}

If $H$ is a Heyting algebra, the \emph{ideal spectrum} $\Spec_{\down}(H)$ 
is the partially ordered set of all prime ideals of $H$ endowed with
the \emph{Zariski topology}, a basis of which consists of the sets 
$\bar I(a) := \{ \mf i\in\Spec_{\down}(H) \mid a\notin\mf i \}$. 
The first part of the following fact goes back to Marshall Stone in
\cite{Stone}.

\begin{Fact} \label{F-spectrum} 
  The map $a \mapsto \bar I(a)$ defines (functorially) an embedding of
  Heyting algebras 
  $$H \ \into\ \mc O(\Spec_{\down}(H)).$$
  Moreover, if $H$ is generated by $g_1,\dots,g_n$, then $\Spec_{\down}(H)$, 
  partially ordered by inclusion, can be turned into a Kripke model of
  \IPL{n} by defining $\val(\mf i) := \{P_i \mid i \in \bar I(g_i)\}$.
\end{Fact}

In the special case of $\Fr{n}$, the Kripke model $\mf K_n$
constructed by Bellissima is naturally isomorphic to the restriction of
this construction to the principal ideal spectrum, as we will show in
the remaining of this section.

\begin{Rem}
There is a bijection $\mf i \mapsto \mf i\cpl$ between the
set of prime ideals $\Spec_{\down}(H)$ and the set of prime filters 
$\Spec^{\up}(H) = \Spec_{\down}(H\dual)$. 
Therefore there is also an embedding 
\begin{align*}
 H \ \into\ & \mc O \big( \overline{\Spec}^\up(H) \big) \\
 a \ \mapsto\ & F(a) \,:=\, \{\mf p \mid \mf p \text{ prime filter}, a\in\mf p\}
\end{align*}
where $\overline{\Spec}^\up(H)$ denotes the space of prime filters
endowed with the co-Zariski topology, a basis of open sets of which is
given by the $F(a)$'s. 
\end{Rem}

A \emph{principal ideal} in a lattice is the decreasing set generated
by an element $x\neq 1$, which we denote by $(x)_\down$. 
A principal ideal is prime if and only if its generator is
$\hcap$-irreducible. Let the
\emph{principal ideal spectrum} $\Spec^0_\down(H)$ be the space of all
principal prime ideals endowed with the (trace of the) Zariski 
topology. The continuous inclusion map $\iota: \Spec^0_\down(H) \to
\Spec_\down(H)$ induces an epimorphism of Heyting algebras 
$\iota^*: \mc O(\Spec_\down(H)) \to \mc O(\Spec^0_\down(H))$, given by
$U \mapsto U \cap \Spec^0_\down(H)$.

\medskip
\textbf{Definition } 
We denote by $\bar I_0$ the map 
$\iota^* \circ \bar I:\ H \to \mc O\big(\Spec^0_\down(H)\big)$.

\begin{Prop} \label{Pr-recoverKripkeb}
For $H$ dense in $\Fco{n}$, we have that 
$$\mc O \big( \Spec^0_\down(H) \big) \ =\ 
  \mc O_\down \big( \Spec^0_\down(H),\subseteq \big),$$ 
that is, the Zariski topology on the principal ideal spectrum is the
topology of decreasing sets.
\end{Prop} 

\PROOF
By Lemma~\ref{L-denseirr}, $H$ has the same $\Hcup$- and
$\Hcap$-irreducibles as $\Fco{n}$. As the latter satisfies the
assumptions of Lemma~\ref{L-dual}, we may use it freely. In fact, the
proposition holds more generally for a Heyting algebra satisfying the 
duality in Lemma~\ref{L-dual}. 

By definition of the Zariski topology on the ideal spectrum, the inclusion
``$\subseteq$'' is clear. Conversely, let $X$ be a decreasing set in
$\Spec^0_\down(\Fr{n})$ with respect to inclusion. We have to show
that $X$ is open in the Zariski topology. In fact, we are proving 
$X = \bigcup \big\{ \bar I_0(a^\hcup) \bigm| (a)_\down \in X \big\}$:

From $a^\hcup \not\hle a$ (Lemma~\ref{L-dual}) it follows that
$a^\hcup \notin (a)_\down$, \textit{i.e.} $(a)_\down \in 
\bar I_0(a^\hcup)$ and thus ``$\subseteq$''. 
Conversely, let $(b)_\down \in \bar I_0(a^\hcup)$ for some 
$(a)_\down \in X$, i.e. $a^\hcup \notin (b)_\down$. This means
$a^\hcup \not\hle b$, whence (by definition of ${}^\hcap$, see
Lemma~\ref{L-dual}) $b \hle a^{\hcup\hcap} = a$. This implies
$(b)_\down \subseteq (a)_\down$, and because $X$ is decreasing,
$(b)_\down \in X$. 
\END

\begin{Rem} \label{Rem-eins}
By Lemma~\ref{L-irr}, in $\Fr{n}$ as well as in $\Fco{n}$, the
$\hcap$-irreducibles are $\Hcap$-irreducible, and they are the
same. Therefore, and with Proposition~\ref{Pr-recoverKripkeb}, 
$$\arraycolsep3pt\begin{array}{@{}ccccc@{}}
 \mc O \big( \Spec^0_\down(\Fr{n}) \big) & = 
  & \mc O_\down \big( \Spec^0_\down(\Fr{n}),\subseteq \big)
  & \isom_{\text{naturally}} &
  \mc O_\down \big(\text{$\textstyle\Hcap$-irreducibles} \text{ of } \Fr{n},\hle \big)\\
  &&&& || \\ 
  \mc O \big( \Spec^0_\down(\Fco{n}) \big) & = 
  & \mc O_\down \big( \Spec^0_\down(\Fco{n}),\subseteq \big)
  & \isom_{\text{naturally}} & 
  \mc O_\down \big(\text{$\textstyle\Hcap$-irreducibles} \text{ of } \Fco{n},\hle \big)
\end{array}$$
In particular, $\Spec^0_\down(\Fco{n})$ and $\Spec^0_\down(\Fr{n})$
are naturally homeomorphic. 
(Note that $\Spec_\down(\Fco{n})$ and $\Spec_\down(\Fr{n})$ are not
homeomorphic.) 
\end{Rem}

\begin{Theorem} \label{Thm-recoverKripke}
  The generic Kripke model $\mf K_n$ is, as a partial order, the principal 
  ideal spectrum of $\Fr{n}$ (equivalently of $\Fco{n}$) with inclusion. The
  valuations are determined   by the images $\bar I_0(g_j)$ of the
  free generators $g_j$ of $\Fr{n}$, namely $P_j \in \val(\mf i)$ for some
  principal ideal $\mf i$ iff $g_j \notin \mf i$, or equivalently, 
  $\mf i \in \bar I_0(g_j)$. 
\end{Theorem}

\PROOF 
By Fact~\ref{F-Bneu}\,(a) and Lemma~\ref{L-irr}, an element of $\Fr{n}$
or $\Fco{n}$ is determined by its $\supp_{\Hcap}$, that is by the
$\Hcap$-irreducibles of the algebra. It follows that the map $\bar
I_0$ is injective, Together with Proposition~\ref{Pr-recoverKripkeb},
we get an embedding 
$\Fr{n} \ \into\ \mc O_\down \big( \Spec^0_\down(\Fr{n}),\subseteq \big)$, 
and the right side is, by the previous remark, naturally isomorphic to 
$\mc O_\down \big(\text{$\textstyle\Hcap$-irreducibles} \text{ of } \Fr{n},\hle \big)$. 
The $\Hcap$-irreducibles are of the form $\Zre w$ for $w \in K_n$ 
by Fact~\ref{F-Bneu}\,(a), and $v \le w \iff \Zre v \hle \Zre{w}$.
This proves the first statement. 

A principal prime ideal of $\Fr{n}$ is generated by some $\Zre w$ for 
$w \in K_n$. As the unique successor of $\Zre w$ is $\Zre w \cup \{w\}$, 
an element $x$ of $\Fr{n}$ is not in the ideal generated by $\Zre w$ 
iff $w \in x$. Thus the free generators $g_j$ are mapped on 
$$\bar I_0(g_j) \ =\
  \{ \mf i \in  \Spec^0_\down(\Fr{n}) \mid g_j \notin \mf i\} 
  \;\ \hat=\;\  \{ w \in K_n \mid w \in g_j\},$$ 
where ``$\hat=$'' stands for the image under the natural isomorphism. 
On the other hand the image of $g_j$ is $\{ w \in K_n \mid P_j \in \val(w)\}$. 
This proves the second part of the theorem.
\END

\begin{Rem} \label{Rem-machin}
Theorem~\ref{Thm-recoverKripke} identifies an element $a$ of
$\Fr{n}$ with $\supp_{\Hcap}(a)$, whereas Bellissima's construction is
more easily understood as identifying it with $\supp_{\Hcup}(a)$, that is
with the underlying embedding 
$$\Fr{n} \ \into\ \mc O \big(\overline{\Spec}_c^\up(\Fr{n})\big) \ \isom\ 
  \mc O_\down \big(\text{$\textstyle\Hcup$-irreducibles},\hle \big),$$
where $\overline{\Spec}_c^\up(\Fr{n})$ is the space of ``completely
prime principal filters''. Algebraically, this space is less natural
than the principal ideal spectrum --- the lack of complete duality
comes from the fact that not all $\hcup$-irreducible elements are
completely $\hcup$-irreducible. 
However, up to the duality of Lemma~\ref{L-dual}, the embedding
is the same as in Theorem~\ref{Thm-recoverKripke}. In the light of 
Remark~\ref{Rem-eins}, one sees that this duality is nothing else 
than the order preserving homeomorphism between $\Spec^0_\down(\Fr{n})$ 
and $\overline{\Spec}_c^\up(\Fr{n})$ mapping a principal prime ideal
on its complement, which is exactly the map $(a)_\down \mapsto (a^\hcup)^\up$. 
\end{Rem}

\section{Some model theory of finitely generated free Heyting algebras} \label{Sec-mt}

The basic model theoretic notions like elementary equivalence
$\equiv$, elementary substructure $\elem$, definability and
interpretability, are explained in any newer model theory
textbook, see for example \cite{Hodges}. ``Definable'' means definable with
parameters, and ``$A$-definable'' with parameters in $A$. 

The theory of Heyting algebras has a model completion (in \cite{GZ},
as a consequence of a result by Pitts \cite{P}), and there are some
results about (un)decidability (see for example \cite{R} and \cite{I}), 
but otherwise little seems to be known about the model theory of
Heyting algebras.

\subsection{First order definition of the Kripke model}

\begin{Theorem} \label{Thm-gendef} 
  Fix free generators $g_1,\dots,g_n$ of $\Fr{n}$. Let $H$ be dense in
  $\Fco{n}$ and containing $g_1,\dots,g_n$. Then the set
  $\{g_1,\dots,g_n\}$ is $\0$-definable in $H$. 
\end{Theorem}

\PROOF 
First we note that the partial order $(K_n,\le)$ of the Kripke model 
$\mf K_n$ is $\0$-definable in $H$: the underlying set can be
identified with the $\Hcup$-irreducibles of $H$ by Lemma~\ref{L-denseirr}. 
It is $\0$-definable as the set of those elements  having a unique
predecessor. They are ordered by the restriction of the partial order
of $H$. According to Remark~\ref{Rem-suppKripke}, this order can be
identified with $(K_n,\le)$. 

In the sequel of the proof, we will simply write $K_n$ for the
definable set of $\Hcup$-irreducibles of $H$. We have then a
$\0$-definable injection $H \to \mf P(K_n)$ that maps an element 
$a$ on its $\{a\}$-definable \emph{support} $\supp_{\Hcup}(a) = 
\{w \in K_n \mid w \hle a\}$. In this proof, ``successor'' and
``predecessor'' are always meant in $(K_n,\le)$.

Clearly, the set of atoms of $\Fr{n}$ is $\0$-definable. For example,
they are exactly the elements whose support is a singleton. The unique
element of the support of an atom $a$ will be called $w_a$. Let $a$ be
an atom, and $\beta$ the valuation of $w_a$. Consider the set of all
elements of $K_n$ of the form $w_{\beta',\{w_a\}}$ for $\beta' \subset
\beta$. It has $2^{|\beta|} - 1$ elements and is $\{a\}$-definable
because it consists of all elements $w \in K_n$ which are successors
of $w_a$ without other predecessors. Therefore for any $k$, the set
$A_k := \big\{ a \bigm| a \text{ atom} \text{ and } |\val(w_a)| = k \big\}$
is $\0$-definable. 

Let $a_i$ be the atom with $w_{a_i} = w_{\{P_i\},\emptyset}$. 
First we remark that the set of atoms $$B_i \,:=\, \big\{ a \bigm| 
a \text{ atom} \text{ and } P_i \in \val(w_a) \big\}$$ 
is $\{a_i\}$-definable, because this is exactly $a_i$ together with 
the set of those atoms $a$ such that the point $w_a$
has two common successors with $w_{a_i}$ without other predecessors,
namely $w_{\0,\{w_a,w_{a_i}\}}$ and $w_{\{P_i\},\{w_a,w_{a_i}\}}$. 
Now $A_1 = \{a_1,\dots,a_n\}$ is a finite $\0$-definable set.
Therefore, to prove the proposition, it is sufficient to show that
$\supp_{\Hcup}(g_i)$ is $\{a_i\}$-definable (uniformly in $i$) .

\medskip\textbf{Claim:} 
$\supp_{\Hcup}(g_i)$ consists of all points $v \in K_n$ satisfying the 
following first order conditions: 

(1) $v$ has a successor in $(K_n,\le)$ that has no other predecessor than $v$; \\ 
(2) either $v \ge w_{a_i}$ or $v$ has two common successors with $w_{a_i}$
   that have no other predecessors. 

\medskip Proof: 
For the inclusion ``$\subseteq$'', note first that any element 
$v \in \supp_{\Hcup}(g_i)$ has the successor $w_{\0,\Erz{v}}$  
that has no other predecessor. Then, if $v \in \supp_{\Hcup}(g_i)$ is
not above $w_{a_i}$, then there are the two elements $w_{\0, \Erz{v,w_{a_i}}}$ 
and $w_{\{P_i\}, \Erz{v,w_{a_i}}}$ satisfying (2).

For the converse inclusion, we first notice that no point with
valuation $\0$ can satisfy condition (1) since the Kripke model is
reduced. If $v \ge w_{a_i}$, then the valuation
of $v$ is either $\0$ or $\{P_i\}$; the former is excluded by (1). 
If $v \not\ge w_{a_i}$ and $P_i \notin \val(v)$, then there is only
one common successor with $w_{a_i}$ without other predecessors, namely 
$w_{\emptyset,\Erz{v,w_i}}$, contradicting (2).
\END

\medskip
\textbf{Definition } 
We call \emph{pre-generators} of $\Fr{n}$ the atoms $a_1,\dots,a_n$
such that $|\val(a_i)| = 1$, and we will fix them for the remaining of
this section. 
(With the notation of the previous proof, the pre-generators are the
elements of $A_1$.)

\medskip
The proof of the theorem shows in particular that in a dense
sub-algebra of $\Fco{n}$ containing $g_i$, the corresponding
pre-generator $a_i$ is interdefinable with $g_i$. On the one hand,
$g_i$ is the unique element having the $\{a_1\}$-definable set
$\supp_{\Hcup}(g_i)$ as its support; on the other hand, $a_i$ is the
unique element in the $\{g_i\}$-definable set $A_1 \cap \supp_{\Hcup}(g_i)$.

\begin{Cor} \label{Cor-intKripke}
  If $H$ is dense in $\Fco{n}$, then the Kripke model $\mf K_n$ is
  interpretable in $H$ with parameters $a_1,\dots,a_n$. 
\end{Cor} 

\PROOF
We have already seen in the proof of Theorem~\ref{Thm-gendef} that
the partial ordering $(K_n,\le)$ is $\0$-definable and that the
support $\supp_{\Hcup}(g_i)$ is $\{a_i\}$-definable. Now the points 
$x \in K_n$ with valuation $\{P_i\mid i \in I\}$ are definable as
those satisfying the formula that expresses 
$\bigwedge_{i=1}^n \big( x \in \supp_{\Hcup}(g_i) \iff i\in I \big)$. 
\END

\begin{Cor}[Grigolia \cite{Grigolia-Aut}] \label{Cor-freegen}
  $\Fr{n}$ has only one set of free generators, and hence
  $\Aut(\Fr{n}) = \mr{Sym}(n)$.
\end{Cor} 

\PROOF 
Any set of free generators of $\Fr{n}$ has size $n$ (because there are
$2^n$ atoms). Assume $\Fr{n}$ is freely generated by $g_1,\dots,g_n$
and $h_1,\dots,h_n$. Then $g_i \mapsto h_i$ extends to an automorphism
of $\Fr{n}$, which has to leave the $\0$-definable set $\{g_1,\dots,g_n\}$ 
invariant. By definition of $\Fr{n}$ as the free algebra, any
permutation of the free generators extends uniquely to an automorphism
of $\Fr{n}$. 
\END

\begin{Cor}  \label{Cor-Aut} 
  If $H$ is dense in $\Fco{n}$, then $\Aut(H) \leqslant \mr{Sym}(n)$.
  If $H$ is in addition setwise invariant under $\Aut(\Fco{n})$, as for
  example $\Fco{n}$ and $\Fkl{n}$, then $\Aut(H) = \mr{Sym}(n)$. 
\end{Cor} 

\PROOF 
The metric on $\Fr{n}$ is invariant under $\Aut(\Fr{n})$, hence every
automorphism is continuous and therefore extends uniquely to the
completion $\Fco{n}$. 
Let $H$ be dense in $\Fco{n}$. As any automorphism of $H$
permutes the $\0$-definable set $\{a_1,\dots,a_n\}$, we get a map
$\Aut(H) \to \mr{Sym}(\{a_1,\dots,a_n\})$. Let $\alpha$ be in the
kernel, i.e. fixing $a_1,\dots,a_n$ pointwise. We have to show that
$\alpha$ is the identity. Now $\alpha$ fixes the Kripke model 
$\mf K_n$ interpreted in $H$ as in \ref{Cor-intKripke}. But every
element of $H$ is interdefinable with a subset of $K_n$, namely its
support. Therefore $\alpha$ has to be the identity. 

Conversely, any automorphism of $\Fco{n}$ restricts to $H$ if $H$ is
invariant, so $\Aut(H) = \mr{Sym}(n)$ in this case. $\Fkl{n}$ is
invariant as being generated by all the $\Hcap$-irreducibles. 
\END

Let $\mr{dcl}_T$ and $\mr{acl}_T$ stand for the definable and model
theoretic algebraic closure in the theory $T$, see e.g.\ \cite{Hodges}.

\begin{Cor} \label{Cor-acl} 
  If $H$ is dense in $\Fco{n}$, then $$\Fr{n} \cap H 
  \;\subseteq\; \mr{dcl}_{\mr{Th}(H)}(a_1,\dots,a_n) 
  \;\subseteq\; \mr{acl}_{\mr{Th}(H)}(\0).$$ 
\end{Cor}

\PROOF
Corollary~\ref{Cor-intKripke} allows us, over the parameters
$a_1,\dots,a_n$, to define the supports of the generators of
$\Fr{n}$. Now every element of $\Fr{n}$ is a term in the generators. 
This implies that the support of every element $x$ in $\Fr{n}$ is
$\{a_1,\dots,a_n\}$-definable (cf. Remark~\ref{Rem-suppminneu}). 
If $x$ is also in $H$, then $x$ is $\{a_1,\dots,a_n\}$-definable as
the unique element having its support. The second inclusion is clear
as the $a_i$ are algebraic over $\0$ (for example as the elements of
the finite $\0$-definable sets of atoms). 
\END

In particular, $\Fr{n} \subseteq \mr{acl}_{\mr{Th}(\Fco{n})}(\0)$.

\begin{Qu}
  Does equality hold?
\end{Qu}

\subsection{Comparing theories}

What can be said about the first order theories of $\Fco{n}, \Fr{n}$
and $\Fkl{n}$? 

As $n$ is coded in the number of atoms, which are first order
definable, we get that $H_n \not\equiv H_m$ if $n\neq m$, $H_n$ is
dense in $\Fco{n}$ and $H_m$ dense in $\Fco{m}$. More precisely,
this proves a difference in the $\forall\exists$-theories. 
Bellissima's Corollary 3.2 in \cite{B} gives a better result, namely 
$(\Fr{n})^{}_\forall \neq (\Fr{m})^{}_\forall$ for $n\neq m$, due to
an ``identity'', i.e. a positive universal formula.
With Proposition~\ref{Pr-elema}, it follows that $(H_n)_\forall
\neq (H_m)_\forall$ for $H_n,H_m$ as above. 

Comparing $\Fco{n}, \Fr{n}$ and $\Fkl{n}$ with the same $n$, 
we have to distinguish the case $n=1$ where $\Fkl{1} = \Fr{1} = \Fco{1}$
from the case $n>1$ where the three algebras are pairwise not isomorphic: 
$\Fco{n}$, has size continuum, whereas $\Fr{n}$ and $\Fkl{n}$ are
countable; $\Fr{n}$ is finitely generated, but $\Fkl{n}$ is not
(Fact~\ref{F-Fsmile}).  

Concerning elementary equivalence and similar concepts, 
Theorem~\ref{Thm-gendef} yields the following results:

\begin{Prop} \label{Pr-minprim}
  $\Fr{n}$ embeds in every model of its theory (``$\Fr{n}$ is an
  algebraic prime model''). No proper dense sub-algebra of $\Fr{n}$ is
  elementarily equivalent to $\Fr{n}$. 
\end{Prop} 


\PROOF 
The generators form a finite $\emptyset$-definable set, thus they
belong to every model. Every element of $\Fr{n}$ is a term in the
generators, and the theory of $\Fr{n}$ knows which terms describe the
same element in $\Fr{n}$ and which not. Therefore $\Fr{n}$ is a
substructure of every model. 

If $H$ is a dense sub-algebra of $\Fr{n}$, then $H$ has the same atoms
as $\Fr{n}$ and interprets the partial ordering $(K_n,\le)$ of the
Kripke model in the same way as $\Fr{n}$. Now $\Fr{n}$ satisfies the
formula saying that there are elements $g_1,\dots,g_n$ such that their
supports in $K_n$ are defined from the atoms in $A_1$ as in the proof
of Theorem~\ref{Thm-gendef}. Thus any elementarily equivalent
dense sub-algebra has to contain the generators. 
\END

With Grigolia's result that no proper sub-algebra of $\Fr{n}$ is
isomorphic to $\Fr{n}$, the first part of the proposition immediately
implies

\begin{Cor} 
No proper sub-algebra of $\Fr{n}$ can be elementarily equivalent to
$\Fr{n}$.
\end{Cor}

We do not know whether $\Fr{n}$ is also an elementary prime model of
its theory (i.e. embeds elementarily into every model of its theory),
and we do not know whether the free Heyting algebra is an elementary
substructure of its completion.

\begin{Prop} 
  If $\Fr{n} \equiv \Fco{n}$, then $\Fr{n} \elem \Fco{n}$.
\end{Prop} 

\PROOF 
If $\Fr{n} \equiv \Fco{n}$, then Theorem~\ref{Thm-gendef} implies 
$$(\Fr{n},g_1,\dots,g_n) \,\equiv\, (\Fco{n},g_{\sigma(1)},\dots,g_{\sigma(n)})$$
for some $\sigma \in \mr{Sym}(n)$. With Corollary~\ref{Cor-Aut} we
then get $$(\Fr{n},g_1,\dots,g_n) \,\equiv\, (\Fco{n},g_1,\dots,g_n).$$ 
Finally the result follows from Corollary~\ref{Cor-acl} and the
interdefinability of $a_i$ and $g_i$. 
\END

\begin{Prop} \label{Pr-elema}
  If $H$ is dense in $\Fco{n}$, then $H \elem_\forall^+ \Fco{n}$,
  which means that both algebras satisfy the same positive universal
  $\mc L_{HA}$-formulae with parameters in $H$. 
\end{Prop} 

\PROOF 
It is clear that if $\Fco{n}$ satisfies a universal formula, then also
$H$. Assume $\Fco{n} \models \exists \bar x\,\phi(\bar x,\bar a)$ where
$\phi$ is a negative quantifier-free formula with parameters $\bar a$
from $H$. Then $\phi$ can be put in the form $\bigwedge_i\bigvee_j
\tau_{ij} (\bar x, \bar a) \neq 1$ for $\mc L_{HA}$-terms
$\tau_{ij}$. Each term defines a continuous function (see
Theorem~\ref{Thm-metric}), and as points are closed, 
$\tau_{ij} (\bar x, \bar a) \neq 1$ defines an open set. Thus
$\phi(\bar x,\bar a)$ defines an open set in $\Fco{n}^l$ where $l$ is
the length of $\bar x$. If this open set is non-empty as the formula
above asserts, then the intersection with the dense subset $H^l$ is
also non-empty. 
\END

In the language of universal algebra, $H \elem_\forall^+ H'$ means
that $H$ satisfies the same identities as $H$ in the language with
constants for all element of $H$. In particular, the proposition
provides a proof of Lemma~4.6 in \cite{B}.%
\footnote{The proof of Lemma 4.6 in \cite{B} uses Lemma 4.5, which
  contains a mistake: The hypothesis must be $a_i \cap H_{\alpha,n} =
  b_i \cap H_{\alpha,n}$. Otherwise (with $w_i$ as in figure~1 p.156 of
  \cite{B}) for $\alpha = 1$, $a_0 = \{w_0,w_2\}$, $b_0 = \{w_0,w_1,w_2\}$ 
  and $p(x) = (x \hto 0) \hto 0$ one gets a counterexample, as 
  $a_0 \cap \mr{Lev}_{1,1} = b_0 \cap \mr{Lev}_{1,1} = \{w_2\}$, but
  $p(a_0) = a_0$ and $p(b_0) = 1$, thus $p(b_0) \cap \mr{Lev}_{1,1} =
  \{w_2,w_3\}$. But the proof of Lemma 4.6 works with this weaker 
  version of Lemma 4.5.}

\begin{Cor} \label{Cor-hcapirr}
  If $H$ is dense in $\Fco{n}$, then every $\hcap$-irreducible
  element of $H$ remains $\hcap$-ir\-re\-du\-cible in $\Fco{n}$.
\end{Cor}

\PROOF 
An element $u \in H$ is $\hcap$-irreducible iff the positive universal
formula $\forall x\, (x \hcup u = u \lor (x \hcup u) \hto u = u)$ holds
in $H$.
\END 

Thus the $\hcap$-irreducible elements of a dense sub-algebra are
exactly the co-principal sets. For $\hcup$-irreducible elements, the
situation is different: the corresponding result of the corollary
holds (see Corollary~\ref{Cor-hcupirr}), but there are more
$\hcup$-irreducibles than just the principle sets, and in general not
all the $\hcup$-irreducibles of $\Fco{n}$ are in a dense sub-algebra.

\begin{Rem} 
  $\Fco{n}$, as the profinite limits of the finite Heyting algebras
  $\Fr[d]{n}$, can be embedded in a pseudo-finite Heyting algebra,
  namely in a nontrivial ultraproduct of the $\Fr[d]{n}$ via 
  $x \mapsto (\pi_d(x))_{d \in \omega} \in \big( \prod_{d \in \omega}
  \Fr[d]{n} \big) /\mc U$. Hence $\Fco{n}$ (and hence every dense
  sub-algebra) satisfies the universal theory of all finite Heyting
  algebras. 

Similarly, if $H$ is a dense sub-algebra of $\Fco{n}$, we can map 
$\Fco{n}$ in an ultrapower of $H$, via 
\begin{align*} 
  \Fco{n} \to H^{\mc U}, \quad
  & x \mapsto \big( \sigma_d^{\min}(\pi_d(x)) \big)_{d\in \omega} \\
  \text{ or via } \qquad 
  & x \mapsto \big( \sigma_d^{\max}(\pi_d(x)) \big)_{d\in \omega}
\end{align*} 
with the sections $\sigma_d^{\min}, \sigma_d^{\max}$ as in
Remark~\ref{Rem-section}. If $\mc U$ is a non-trivial ultrafilter on
$\omega$, then these are $\{0,\hcap,\hcup,\hle\}$-embeddings and
$\{0,1,\hcap,\hto,\hle\}$-embeddings respectively. Thus $H$ has 
the same universal theory as $\Fco{n}$ in any of the two languages: 
$\{0,\hcap,\hcup,\hle\}$ and $\{0,1,\hcap,\hto,\hle\}$. 
\end{Rem}

\section{Further remarks and open problems} \label{Sec-letzt}

\subsection{Open problems}

\begin{Problem} 
  Is it possible to characterise the subsets of $K_n$ that are in
  $\Fr{n}$? in $\Fkl{n}$?
\end{Problem}

By Fact~\ref{F-spectrum}, every Heyting algebra embeds into a
topological Heyting algebra. Therefore, the universal theory of all
topological Heyting algebras equals $(T_{HA})_\forall$. In particular, 
on the quantifier-free level one can compute in the theory of Heyting
algebras as if one were in an arbitrary topological space.

\begin{Problem} \label{Prob-top}
Does $T_{HA}$ equal the theory of all topological Heyting algebras? 
\textit{I.e.}\ does any $\mc L_{HA}$-sentence which holds in all
lattices of open sets of topologies hold in all Heyting algebras? 
\end{Problem}

\begin{Problem} 
 Is $\Fr{n} \elem \Fco{n}$? 
 Does $\Fr{n}$ eliminate quantifiers in a reasonable language?
\end{Problem}

\subsection{The $\hcup$-irreducible elements} \label{S-hcupirr}

Fact~\ref{F-Bneu}, Lemmas~\ref{L-irr}, \ref{L-denseirr} and 
Corollary~\ref{Cor-hcapirr} completely determine the $\Hcup$-, $\Hcap$- 
and $\hcap$-irreducible elements of $\Fco{n}$ and its dense sub-algebras. 
Now we are going to characterise the $\hcup$-irreducible elements. 
For $n=1$, the principal sets and $1$ are the only $\hcup$-irreducibles; 
for $n>1$, there are more infinite $\hcup$-irreducibles.

\begin{Prop} \label{Pr-hcupirr}
  If $H$ is dense in $\Fco{n}$, then $X \in H$ is $\hcup$-irreducible
  iff for all (incomparable) $w_0,w_1 \in X$ there exists an element
  $w \in X$ with $w \geqslant w_0$ and $w \geqslant w_1$, i.e. $X$ as
  a subset of $(K_n,\le)$ is upward filtering. 
\end{Prop}

\PROOF
If $X$ is a proper union of $X_0,X_1 \in H$, choose $w_i \in X_i
\setminus X_{1-i}$. Then they are incomparable and have no common
larger element $w$ in $X$. 
Conversely, let $w_0,w_1 \in X$ and define $X_i := X \cap \Zre{w_i}$.
Then $X_i$ is a proper subset of $X$, $X_i \in H$ because the
co-principle sets are in $H$ by Lemma~\ref{L-denseirr}, and 
$X_0 \cup X_1 = X \cap \big( \Zre{w_0} \cup \Zre{w_1} \big) = 
X \setminus \big( \{w_0\}^\up \cap \{w_1\}^\up \big)$. If $X$ is
$\hcup$-irreducible, then $X_0 \cup X_1 \neq X$, and there is 
$w \in X \cap \{w_0\}^\up \cap \{w_1\}^\up$. 
\END

\begin{Cor} \label{Cor-hcupirr}
  If $H$ is dense in $\Fco{n}$, then a $\hcup$-irreducible element of
  $H$ remains $\hcup$-irreducible in $\Fco{n}$. 
\end{Cor}

\begin{Prop} 
For $n>1$, there are continuum many $\hcup$-irreducibles in $\Fco{n}$.
\end{Prop} 

\PROOF 
There exists an infinite antichain $(z_i)_{i\in\omega}$ in $K_n$
(Fact~\ref{F-anti}), and for any proper subset $I$ of $\omega$, the
set $Z_I := \bigcap_{i \in I} \Zre{z_i}$ is $\hcup$-irreducible by
Proposition~\ref{Pr-hcupirr}: 
for $w_0,w_1 \in Z_I$ and $j \notin I$, there is a common larger
element $w_{\emptyset, \Erz{w_0,w_1,z_j}}$. 
\END

In particular, not every $\hcup$-irreducible element of $\Fco{n}$ is
in $\Fr{n}$. Also it follows from this proof that all co-principal sets
are $\hcup$-irreducible (Theorem~3.1 in \cite{B}), because any element
of $K_n$ is part of a two-element antichain.

It is easy to check that the element $\bigcup_{i\in\omega} \Erz{z_i}$ 
of $\Fco{n}$ is not a finite union of $\hcup$-irreducible elements. 
In contrast to this, Urquhart (Theorem~3 in \cite{Urquhart}) has shown
that every element of $\Fr{n}$ is a finite union of $\hcup$-irreducible 
elements.

\begin{Qu} 
Is $\Fco{n}$, or more generally any dense sub-algebra, generated by
its $\hcup$-irreducible elements? 
\end{Qu}

\begin{Prop} 
Any intersection of some of the free generators of $\Fr{n}$ is
$\hcup$-ir\-re\-du\-cible in both, $\Fr{n}$ and $\Fco{n}$. 
\end{Prop} 

\PROOF
Consider $[\![P_1]\!] \cap\dots\cap [\![P_k]\!]$, i.e. all
points whose valuation includes $P_1,\dots,P_k$. For any two such
points $w_0,w_1$, either they are comparable and the larger one is a
common larger element, or they are incomparable, and then 
$$w \,=\, w_{\{P_1,\dots,P_k\},\Erz{w_0,w_1}} \in\; 
  [\![P_1]\!] \cap\dots\cap [\![P_k]\!]$$ is a common larger element. 
\END

It follows that if $n \ge 2$, then any intersection of at most $n-1$
of the generators is an example of a $\hcup$-irreducible element that
is neither $\Hcup$- nor $\hcap$-irreducible. The intersection of all 
generators is the atom $\{w_{\{P_1,\dots,P_n\},\emptyset}\}$.

\subsection{Approximations of $\protect\Fkl{n}$} \label{S-approx}

Let, as in the proof of Theorem~4.4 in \cite{B}, $B_{n,d}$ be the
sub-algebra of $\Fkl{n}$ generated by all principal sets $\Erz w$ with
$w$ of foundation rank $\le d$; and let $C_{n,d}$ be the sub-algebra
of $\Fkl{n}$ generated by all co-principal sets $\Zre w$ with $w$ of
foundation rank $\le d$. Recall from the proof of
Lemma~\ref{L-denseirr} that  
$$\Zre w = \Erz w \hto \Big( \bighcup_{v<w} \Erz v \Big) \quad\text{and}\quad
  \Erz w = \bighcap \Big\{ \Zre{v} \bigm| v \text{ minimal}\notin \Erz w \Big\},$$
hence we get 
$$C_{n,d} \ \subseteq\ B_{n,d} \ \subseteq\ C_{n,d+1} 
  \ \subseteq\ \cdots \bigcup_{d\in\omega} B_{n,d} 
  \ =\ \bigcup_{d\in\omega} C_{n,d} \ = \ \Fkl{n}.$$
If $n>1$, the inclusions are all strict: 
Bellissima has shown that $B_{n,d}$ can't separate points 
$w_{\beta,Y},w_{\beta',Y} \in K_n^{d+1} \setminus K_n^d$ with
$\beta\neq\beta'$, but the set $\Zre{w_{\beta,Y}}$ in $C_{n,d+1}$
does.  

The proof for the second sort of inclusion is similar, but even
easier: For $w = w_{\{P_i\},Y} \in K_n$ of foundation rank $d$, the
set $\Erz w \in B_{n,d}$ separates $w$ from $w' := w_{\0,\Erz{w}}$. 
On the other hand, $C_{n,d}$ can't separate between $w$ and $w'$: 
This is clear for the generators and clearly preserved under $\hcup$
and $\hcap$, and it is not hard to see that it is also preserved under
$\hto$.

\subsection{Cantor--Bendixson analysis}

$\Fr{1}$ only consists of finite elements and $1$. Thus
Lemma~\ref{L-isolated} implies that the metric topology is the
one-point compactification of a countable discrete set; all points are
isolated, i.e. have Cantor--Bendixson rank $0$, except the maximum 
with rank $1$.

\begin{Prop}
  For $n>1$, $\Fco{n}$ has infinitely many points of rank $1$. They
  are the maximal elements of sub-lattices that look similar to $\Fr{1}$. 
  The elements of higher rank form a perfect subset. 
\end{Prop} 

\PROOF[ (sketchy)]
Let $U^d_{a} := \{x \in \Fco{n} \mid x \cap K_n^d = a\}$ be a
basic open set. By an \emph{extension} of $a$, we mean an $x \in 
U^d_{a}$, and by a \emph{$k$-extension}, we mean an extension $x$ by
adding $k$ new points of the Kripke model. One can check that there
are only the following three possibilities: 

(A) For some $k$, there is no $k$-extension of $a$. Then $U^d_{a}$ 
consists of finitely many finite sets. 

(B) For each $k$, there are exactly two $k$-extensions of $a$.  Then
the extensions of $a$ form a copy of $K_1$, i.e. the Kripke model
for one free generator. Therefore $U^d_{a}$ contains infinitely many
finite sets and exactly one infinite set, which thus is an element of
Cantor--Bendixson rank $1$. 

(C) There at least three $1$-extensions of $a$. Then the number of
$k$-extensions of $a$ increases with $k$. In this case, there are
infinitely many elements in $U^d_{a}$ of rank $>1$ (for each
big enough $l>1$, take a point $w \in K_n^l \setminus K_n^{l-1}$
appearing in some extension of $a$ and then consider the maximal
extension of $a$ omitting this point). 
\END 

\subsection{Order Topologies}

One might wonder how the metric topologies on $\Fr{n}$ and $\Fco{n}$
relate to topologies induced by the partial order $\hle$. There are at
least three topologies that one might consider on a partially ordered
set $(X,\le)$: 
\begin{itemize}
\item
  the topology $\mc O_\down(X,\le)$ of decreasing sets; 
\item
  the topology $\mc O^\up(X,\le)$ of increasing sets; 
\item
  the ``order topology'' $\mc O(X,\le)$ generated by the
  generalised open intervals $(a,b) := \{ x \in X \mid a < x < b\}$ 
  as a sub-basis, where $a = -\infty$ and $b=\infty$ are allowed. 
\end{itemize}

\begin{Prop} 
\begin{enumerate}[{\bf(a)}]
\item
  The trace on $\Fr{n}$ of the increasing and decreasing topology on
  $\Fco{n}$ is the corresponding topology on $\Fr{n}$. 
\item
  For both $\Fr{n}$ and $\Fco{n}$, the order topology contains
  $\mc O_\down$. Neither the order topology nor the metric topology
  contains $\mc O^\up$. 
\item
  For $n=1$, the metric topology on $\Fr{1}$ equals the order
  topology. For $n>1$ and $\Fr{n}$ as well as $\Fco{n}$, the metric
  topology is incomparable with the order topology and does not
  contain $\mc O_\down$. 
\end{enumerate} 
\end{Prop} 

\PROOF 
(a) The intersection of an in-/de-creasing set of $\Fco{n}$ with
$\Fr{n}$ is an in-/de-creasing set of $\Fco{n}$. 

(b) First statement: every proper decreasing set $a$ has successors 
$a \cup \{w\}$ where $w$ is an element of the Kripke model of minimal
foundation rank among those not in $a$. Either there is a unique such
successor, then $a = \big( -\infty,a \cup \{w\} \big)$, or there are
at least two such points $w_1,w_2$ and then $a = \big( -\infty,a \cup
\{w_1\} \big) \cap \big( -\infty,a \cup \{w_2\} \big)$. 

Second statement: $1$ is an isolated point in $\mc O^\up$, but neither
in the order topology, as it does not have predecessors, nor in the
metric topology. 

(c) First statement: First we show that the finite elements are
isolated in the order topology. Let $a$ be finite. If $a=0$, then 
$\{a\} = (-\infty,c)$ for some successor $c$ of $0$. If $a\neq 0$ has
two distinct predecessors $b_1,b_2$, then choose a successor $c$ of
$a$ and then $\{a\} = (b_1,c) \cap (b_2,c)$. Otherwise $a\neq 0$ has a
unique predecessor $b$, that is $a$ is principal. But then $a$ is not
co-principal, hence has two distinct successors $c_1,c_2$ and $\{a\} =
(b,c_1) \cap (b,c_2)$. As remarked in the proof of (a), $1$ is not
isolated in the order topology, therefore the order topology on
$\Fr{1}$ is the same as the metric topology: discrete on the finite 
elements and $1$ is a compactifying point. 

Second statement: 
Consider a metric neighbourhood $\pi_i^{-1}(x)$ containing a
co-prin\-ci\-pal set $\Zre w$ with $w$ of foundation rank less than $i$. 
A neighbourhood of $\Zre w$ in the order topology contains one of the
form $\big( a_1,(\Zre w)^+ \big) \cap\dots\cap \big( a_k,(\Zre w)^+\big)$, 
and because $\Zre w$ does not have predecessors, such a neighbourhood
always contains elements $b$ with $\pi_i(b) = x \cup \{w\}$. 

For the converse and the third statement, consider the decreasing set
generated by $[\![P_1]\!]$. It is also open in the order topology as
it equals $\big( {-\infty}, [\![P_1]\!] \cup \{w_{\{P_2\},\emptyset}\} \big) 
\cap \big( {-\infty}, [\![P_1]\!] \cup \{w_{\emptyset,\emptyset}\} \big)$. 
But it is not open in the metric topology, because its element
$[\![P_1]\!]$ does not contain a metric open neighbourhood: for every finite 
part $a_i := \pi_i([\![P_1]\!])$ the set $a_i\cup\{w_{\emptyset,a_i}\}$ 
is not in $[\![P_1]\!]$, but in $\pi_i^{-1}(a_i)$. (Note that the
existence of $w_{\emptyset,a_i}$ needs $n>1$). 
\END

\begin{Qu} 
  Is the order topology on $\Fr{n}$ the trace of the order topology on
  $\Fco{n}$? 
\end{Qu}



\begin{thebibliography}{BGMM} 
\bibitem[Be]{B}
  Fabio Bellissima, Finitely generated free Heyting algebras, 
  \textit{JSL} \tb{51} (1986) no.\,1, pp.\,152--165. 
\bibitem[BGMM]{BGMM} 
  Guram Bezhanishvili, Mai Gehrke, Ray Mines and Patrick J. Morandi,
  Profinite Completions and Canonical Extensions of Heyting Algebras,
  \textit{Order} \tb{23} (2006) no.\,2--3, pp.\,143--161. 
\bibitem[BMV]{BMV} 
  Patrick Blackburn, Maarten de Rijke, Yde Venema 
  \textit{Modal logic}. 
  Cambridge University Press, Cambridge, 2001.
\bibitem[Bz]{Nick}
  Nick Bezhanishvili
  \textit{Lattices of intermediate and cylindric modal logics}.
  Doctoral thesis, Universiteit van Amsterdam, 2006. 
\bibitem[Bz$^2$]{Bz2} 
  Guram Bezhanishvili and Nick Bezhanishvili, Profinite Heyting
  Algebras, \textit{Order} \textbf{25} (2008) no.\,3, pp.\,211-223. 
\bibitem[Da]{D-preprint}
  Luck Darni\`ere, Model-completion of scaled lattices,
  LAREMA-Preprint No.191, Universit\'e d'Angers, mai 2004.
\bibitem[DJ2]{DJ2}
  Luck Darni\`ere, Markus Junker, Codimension and pseudometric in
  (dual) Heyting algebras, preprint 2008,
  \texttt{http://arxiv.org/archive/math}. 
\bibitem[Fi]{Fitting} 
  Melvin Fitting, \textit{Intuitionistic logic model theory and forcing}.
  North Holland, Amsterdam 1969. 
\bibitem[Gh]{G} 
  Silvio Ghilardi, Free Heyting algebras as bi-Heyting algebras, 
  \textit{C.\ R.\ Math.\ Rep.\ Acad.\ Sci.\ Canada} \textbf{14}
  (1992) no.\,6, pp.\,240--244,  
\bibitem[GhZ]{GZ} 
  Silvio Ghilardi, Marek Zawadowski, 
  Model completions and $r$-Heyting categories, 
  \textit{APAL} \tb{88} (1997), pp.\,27--46. 
\bibitem[Gr1]{Grigolia-buch} 
  Revaz Grigolia, \textit{Free algebras of nonclassical logics}.
  Tbilisi 1987. 
\bibitem[Gr2]{Grigolia-free} 
  Revaz Grigolia, Free and projective Heyting and monadic Heyting
  algebras, pp.\,33--52 in \textit{Non-classical Logics and their
    applications to fuzzy subsets} (ed. U.\,H\"ohle and
  E.\,P.\,Klement), Kluwer Acad. Publ. 1995. 
\bibitem[Gr3]{Grigolia-Aut} 
  Revaz Grigolia Free Heyting algebras and their automorphism groups,
  \textit{Proceedings of  Institute of Cybernetics}, \textbf{2} (2002)
  no.\,1-2 . 
\bibitem[Ho]{Hodges} 
  Wilfrid Hodges \textit{Model theory}. Encyclopedia of Mathematics
  and its Applications \tb{42}, 
  Cambridge University Press, Cambridge, 1993.
\bibitem[Id]{I}
  Pawe\l{}\ Idziak, Elementary theory of free Heyting algebras.
  \textit{Rep.\ Math.\ Logic} \textbf{23} (1989), pp.\,71--73 (1990). 
\bibitem[McT]{McTarski}
  J. McKinsey, Alfred Tarski, On closed elements in closure algebras,
  \textit{Annals of Mathematics} \textbf{47} (1946) no.\,1, pp.\,122--162.
\bibitem[Pi]{P} 
  Andrew Pitts, On an interpretation of second order quantification 
  in first order intuitionistic propositional logic, 
  \textit{JSL} \tb{57} (1992) no.\,1, pp.\,33--52. 
\bibitem[Ry]{R}
  V.V. Rybakov, The elementary theories of free topo-{B}oolean and
  pseudo-{B}oolean algebras, \textit{Mat. Zametki} \textbf{37}
  (1985) no.\,6, pp.\,797--802.
\bibitem[St]{Stone} 
  M. H. Stone, Topological representations of Distributive Lattices
  and Brouwerian Logics, \textit{\v{C}asopis pro \v{p}estov\'an\'i
  matematikyv a fysiky} \tb{67} (1937), pp.\,1--25. 
\bibitem[Ur]{Urquhart}
  A. Urquhart, Free Heyting Algebras, \textit{Algebra Univ.}
  \textbf{3} (1973), pp.\,94--97. 
\end{thebibliography}
\end{document}